\documentclass[pdflatex,sn-mathphys-num]{sn-jnl}%
\usepackage{graphicx}%
\usepackage{multirow}%
\usepackage{amsmath,amssymb,amsfonts}%
\usepackage{amsthm}%
\usepackage{mathrsfs}%
\usepackage[title]{appendix}%
\usepackage{xcolor}%
\usepackage{textcomp}%
\usepackage{manyfoot}%
\usepackage{booktabs}%
\usepackage{algorithm}%
\usepackage{algorithmicx}%
\usepackage{algpseudocode}%
\usepackage{listings}%

\newtheorem{remark}{Remark}%

\usepackage{bm}
\usepackage{bbm}
\usepackage{color}
\usepackage{tabularray}
\usepackage{caption}
\usepackage{diagbox}
\usepackage{csvsimple}

\newcommand{\bff}{\mathbf{f}}

\newcommand{\bfh}{\mathbf{h}}

\newcommand{\bfq}{\mathbf{q}}
\newcommand{\bfr}{\mathbf{r}}

\newcommand{\bfv}{\mathbf{v}}
\newcommand{\bfw}{\mathbf{w}}

\newcommand{\bfA}{\mathbf{A}}

\newcommand{\bfD}{\mathbf{D}}

\newcommand{\bfH}{\mathbf{H}}
\newcommand{\bfI}{\mathbf{I}}

\newcommand{\bfM}{\mathbf{M}}

\newcommand{\bfR}{\mathbf{R}}
\newcommand{\bfS}{\mathbf{S}}

\newcommand{\bfW}{\mathbf{W}}

\newcommand{\barE}{\bar{E}}

\newcommand{\lambdaS}{{\bm \Lambda_{\bf S}}}
\newcommand{\lambdaD}{{\bm \Lambda_{\bf D}}}
\newcommand{\bGamma}{{\bm \Gamma}}

\newcommand{\real}{\mathbb{R}}

\newcommand{\qn}[1][n]{\bfq^{[#1]}}

\DeclareMathOperator{\diag}{diag}

\begin{document}


\title[A Unified Framework for Sponge-Layer Relaxation Methods and Damping Operators for 
Conservation Laws: Application to the Piston Problem of Gas Dynamics]{A Unified Framework for Sponge-Layer Relaxation Methods and Damping Operators for 
Conservation Laws: Application to the Piston Problem of Gas Dynamics}

\author*[1]{\fnm{Carlos } \sur{Muñoz-Moncayo }}\email{carlos.munozmoncayo@kaust.edu.sa}

\affil*[1]{\orgdiv{Computer, Electrical and Mathematical Sciences and Engineering Division}, \orgname{King Abdullah University of Science and Technology (KAUST)}, \orgaddress{\city{Thuwal}, \postcode{23955-6900}, \state{Makkah}, \country{Saudi Arabia}}}

\abstract{

This work addresses the imposition of outflow boundary conditions for one-dimensional 
conservation laws.
While a highly accurate numerical solution can be obtained in the interior of the domain, 
boundary discretization can lead to unphysical reflections.
 We investigate and implement 
some classes of relaxation methods and
far-field operators
to deal with this problem without significantly increasing the size of the computational domain.
We formulate these methods within a framework that allows to reveal relationships among them,
and to propose novel extensions.
 In particular, we introduce a simple and 
 robust relaxation method with a matrix-valued weight function
 that selectively absorbs outgoing waves.
As a challenging model problem, we consider the Lagrangian formulation of the Euler equations for a polytropic gas with inflow boundary conditions determined by an oscillating piston.

}


\maketitle

\section{Introduction}
\label{sec: introduction}

Hyperbolic systems of partial differential equations (PDEs) are characterized
by finite propagation speed of information in their solutions.
At a fixed time, and with given bounded and compactly supported initial data,
we expect the relevant information of a phenomenon described by a system of conservation 
laws to be concentrated inside a bounded subset of the domain.
Nevertheless, when the physical domain is infinite, increasingly extensive computational domains 
might be necessary in numerical simulations,
which can become prohibitively expensive.
Even though this has been partially alleviated by the fast development of computer hardware,
to accurately capture the long-time behavior of solutions still poses a significant challenge.
As a consequence, a variety of techniques have been proposed to reduce spurious reflections at
artificial boundaries (see, e.g.,
\cite{giles_non-reflecting_nodate,hu_absorbing_2004} and references therein).

\begin{figure}
    \centering
    \includegraphics[width=0.7\textwidth]
    {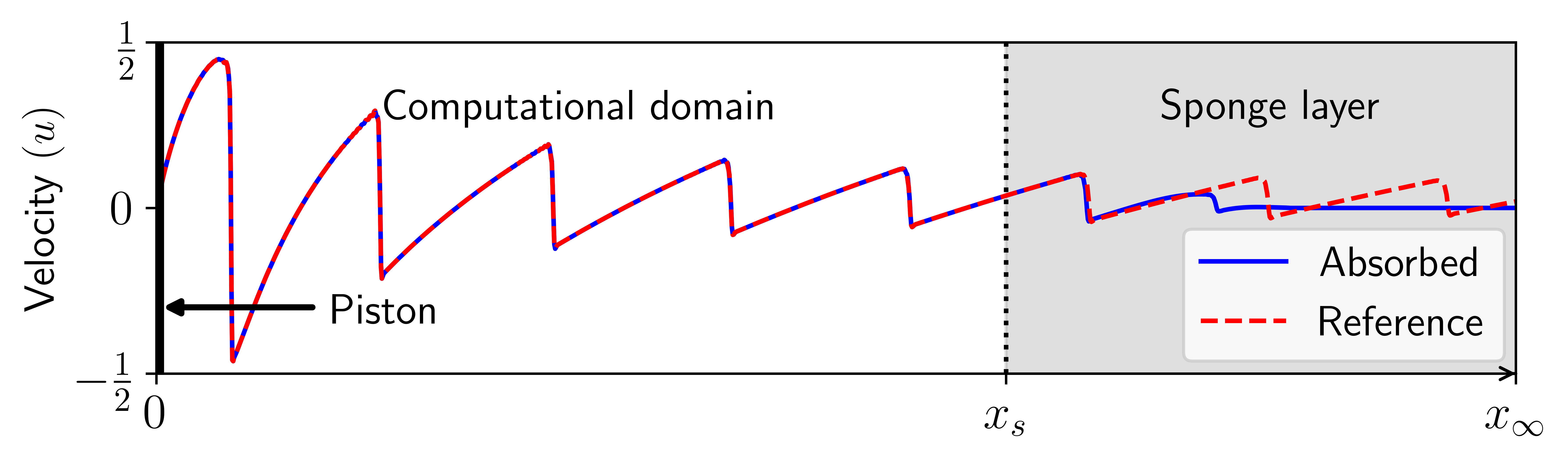}
    \caption{Outgoing waves are absorbed in the sponge layer through the application of artificial boundary conditions.}
    \label{fig:Setup Boundary treatment}
\end{figure}

In 
this manuscript, we focus on 
the so-called
\emph{relaxation} and \emph{far-field operators} 
methods, briefly outlined below.
These are families of 
absorbing boundary conditions (ABCs) \footnote{%
Although the term “absorbing boundary condition” (ABC) is traditionally
reserved for conditions imposed directly at the boundary, we adopt the broader
engineering usage and refer to the zonal approaches considered in this work as
ABCs.}
that could, in principle,
be applied to any system
of hyperbolic PDEs.
Both of them rely on the use of a 
region of the 
domain where the numerical solution
is damped or forced towards a given function.
To avoid ambiguity, we will refer to this region as the sponge layer 
throughout this work (see Figure \ref{fig:Setup Boundary treatment}).

Larsen and Dancy \cite{larsen_open_1983}  introduced the relaxation method 
to apply ABCs to their discretization of Peregrine's
Boussinesq type equations.
The method was first presented as a postprocessing technique where the numerical solution
is divided by an increasing set of positive numbers inside the sponge layer (often referred to
as \emph{relaxation zone} in the literature) at each step
of the time integration procedure.
A Fourier mode-based reflection analysis of a linearization of Peregrine's equations suggested that the
absorbing properties of this technique improve drastically as the ratio between the length of the sponge layer and the wavelength increases.
Due to the versatility and efficiency of this method, it has been widely used  for both absorption and generation of waves in water wave modeling,
where validation tests often consist of solitary waves with a known characteristic length  that
can be contained inside the sponge layer \cite{madsen_boussinesq-type_2003, escalante_efficient_2019, jacobsen_wave_2012}.
Furthermore, several works that formalize this technique and analyze its non-reflecting properties have been proposed, e.g. 
\cite{mayer_fractional_1998,engsig-karup_nodal_2006, engsig-karup_unstructured_nodate, hu_numerical_2015}.

Another popular approach for the development of ABCs is the
modification of the system of governing equations, such
that damping of outgoing waves is a characteristic behavior
of the analytic solution.
In their seminal work \cite{israeli_approximation_1981}, Israeli and Orszag
discuss the addition of reaction and diffusion terms with variable strength
inside the sponge layer. 
They further explore the idea of absorbing waves in just one direction for the wave
equation by modifying its second-order scalar and first-order system forms.
However, it was later shown by Karni \cite{karni_far-field_1996} that such modifications
unavoidably generate reflections that propagate into the computational domain.
Instead, Karni  proposes the use of far-field slowing-down and damping operators
that, by preserving the eigenstructure of the flux Jacobian, absorb waves
in a prescribed direction for systems of hyperbolic PDEs.
A generalization of Karni's slowing-down and damping operators for linear systems
was introduced by Appel\"o and Colonius \cite{appelo_high-order_2009}, 
where high-order viscosity terms are added in the sponge layer.

As a benchmarking application for the ABCs discussed in this work,
we consider the
one-dimensional
Euler
equations in their Lagrangian formulation.
These equations are used to model the 
piston problem of gas dynamics, where the solution is known to 
consist of quasi-periodic nonlinear plane waves with
strength proportional to the amplitude of the 
source 
\cite{webster_finiteamplitude_1977}.
In particular, we study the strongly nonlinear case,
where shock formation occurs in the near-field and
recurrent interaction between the first and second shocks
of the wave train takes place in the far-field, as demonstrated by Inoue and Yano \cite{inoue_propagation_1993}.
A central scheme discretization for this problem and its multi-fluid counterpart 
in Lagrangian coordinates
 was put forward by Fazio and Russo \cite{fazio_central_2010},  where 
the prescription 
of second-order accurate inflow boundary conditions at the piston,
used herein, is described.

The specific contributions of the present work are:
\begin{itemize}
    \item the formulation of a unified framework connecting relaxation methods and
    far-field (damping-operator) approaches for general hyperbolic systems of
    conservation laws (see Sections~\ref{sec:relaxation and far-field} and \ref{sec: Strang});
    \item  a generalization of the relaxation method to improve
    its absorbing properties without compromising robustness or significantly increasing
    its computational cost (see Section \ref{sec: relaxation matrix valued});
    \item a new damping operator based on nonlinear characteristic variables (see Section \ref{sec:NDO});
    \item a quantitative comparison of wave reflection for several state-of-the-art 
ABCs applied to the Lagrangian formulation of the Euler equations, with 
emphasis on the piston problem of gas dynamics (see Section \ref{sec:Results});
\end{itemize}

This manuscript is organized as follows: Section \ref{sec:Boundary treatment}  reviews existing sponge layer-based ABCs.
Section \ref{sec: Extensions} introduces extensions of these techniques, including improved relaxation methods and new damping operators.
Section \ref{sec:Euler equations} recalls the compressible Euler equations in Lagrangian form and the piston problem.
Section \ref{sec:Results} provides a performance comparison of the different ABCs.

The code used for the development of this work is available in the reproducibility repository 
\cite{Munoz_repository_2024}.

\section{Sponge layer-based ABCs}
\label{sec:Boundary treatment}

Unless stated otherwise, hereinafter we consider the governing equations
to be a general hyperbolic system of PDEs of the form
\begin{equation}
    \label{eq:conservation law}
    \partial_t \bfq + \partial_x \bff(\bfq) = {\bf 0}.
\end{equation}
Here $\bfq \in L^1([0,T]\times [0,+\infty[, \real^N)$, and the
Jacobian matrix of $\bff$ can be diagonalized as $\bff'(\bfq)=\bfR\boldsymbol{\Lambda}{\bfR}^{-1}$,
where the columns of $\bfR$ are the right eigenvectors
$\{\bfr_i(\bfq)\}_{i=1}^N$ of $\bff'(\bfq)$, and
$\boldsymbol{\Lambda}$ is a diagonal matrix containing 
its real eigenvalues $\{\lambda_i(\bfq)\}_{i=1}^N$.
The linear case, where $\bff(\bfq)=\bfA \bfq$
for
some constant matrix $\bfA\in\real^{N\times N}$,
\begin{equation}
    \label{eq:conservation law linear}
    \partial_t \bfq + \bfA\partial_x\bfq = {\bf 0},
\end{equation}
will be of relevance for some of the techniques discussed below.
Throughout this work, it is assumed that an appropriate discretization technique is available for these systems,
and a particular approach is briefly discussed in Section \ref{sec:Numerics}.

The computational domain and sponge layer will be the intervals $[0,x_s[$ and $[x_s,x_\infty]$ respectively, and
we suppose there is a far-field
steady-state $\overline{\bfq}\in \real^N$ such that
$\bfq(0,x)=\overline{\bfq}$ for $x>x_s$.
Our goal is to impose ABCs in the sponge layer so that perturbations of $\overline{\bfq}$ in the solution are absorbed before they reach the right boundary $x_\infty$, as depicted in Figure \ref{fig:Setup Boundary treatment}.
Regarding the left boundary, we assume the existence of a function $\bfq_{BC}(\cdot)$
such that well-posedness is guaranteed with the inflow boundary
condition $\bfq(t,0)=\bfq_{BC}(t)$.

\subsection{Slowing-down and damping far-field operators (SDO)}
\label{sec: SDO}
In the work \cite{karni_far-field_1996}, Karni proposes
to solve the linear system
\begin{equation}
        \label{eq: linear SDO}
            \partial_t\bfq + \bfS(x) \partial_x \bfq  = -\bfD(x) \left(\bfq-\overline{\bfq}\right),
\end{equation}
instead of the original system \eqref{eq:conservation law linear}
to impose ABCs.
Here, $\bfS(x)$ and $\bfD(x)$
are referred to as slowing-down and damping far-field operators respectively, and are
given by
\begin{align}
\label{eq: definition far-field operators}
            \bfS(x) = \bfR \lambdaS(x) \bfR^{-1}, \quad \bfD(x) = \bfR \lambdaD(x) \bfR^{-1},
\end{align}
where $\lambdaS(x)$ and $\lambdaD(x)$ are diagonal matrices
with entries
$\{s_{i}(x)\lambda_{i}\}_{i=1}^N$ and $\{d_{i}(x)s_{i}(x)|\lambda_{i}|\}_{i=1}^N$.
The idea behind this approach is to define the operators $\bfS$ and $\bfD$ such
that the solutions of \eqref{eq: linear SDO} coincide
with those of \eqref{eq:conservation law linear}
in the computational domain, while \emph{outgoing waves} are slowed down and damped
in the sponge layer.
In \cite{karni_far-field_1996}, this is achieved by assuming
each
$d_i$ ($s_i$) to be a piecewise constant function
that is $0$ ($1$) in the computational domain and
increases (decreases) to a prescribed value in the sponge layer.
Herein, we motivate this approach at the continuous level,
which will be of further use
in Section \ref{sec: reflection analysis}.

Let $\Gamma_d, \Gamma_s\in \mathcal{C(\real)}$ be two smooth functions that strictly decrease from $1$ to $0$ in the sponge layer,
and are constant everywhere else.
Next, for each $i\in\{1,\ldots,N\}$, we define
\begin{align}
     \label{eq: definition d and s SDO}
    d_i(x)= \sigma_i (1-\Gamma_d(x)), \quad
    s_i(x)= \delta_i \Gamma_s(x),
\end{align}
where $\sigma_i\geq 0$ is a dimensional parameter (with units of inverse time) that
indicates the maximum damping rate that will be attained in the sponge layer,
and $\delta_i\in \{0,1\}$ is an indicator variable that allows to switch off the slowing-down effect
for the $i$-th characteristic field if desired.
By writing System \eqref{eq: linear SDO} in terms of the characteristic variables $\bfw=\bfR^{-1}\bfq$,
we get $N$ initial value problems of the form
\begin{align}
    \label{eq:each characteristic}
    \partial_t w_i + \lambda_i s_i(x) \partial_x w_i = -d_i(x)s_i(x)|\lambda_i| w_i,
    \quad w_i(0,x) = \mathring{w}_i(x) = \left(\mathbf{R}^{-1} \bfq(0,x)\right)_i.
\end{align}
Let us track the solution of \eqref{eq:each characteristic}
accross the characteristic curve $X(t)$ such that 
\begin{align}
     \label{eq:implicit characteristic}
    X'(t) =\lambda_i s_i(X(t)), \quad X(0)=x_0.
\end{align}
This gives us the separable ODE
\begin{equation}
    \label{eq:ODE characteristic}
    \frac{d}{dt}w_i(t,X(t)) = -|\lambda_i| d_i(X(t))s_i(X(t))w_i(t,X(t)),
\end{equation}
which can be solved explicitly to yield
\begin{equation}
    w_i(t,X(t)) = \mathring{w}_i(x_0)  \exp\left(-|\lambda_i| \int_0^t d_i(X(r))s_i(X(r)) \, dr\right).
\end{equation}
For the sake of clarity, let us assume that $s_i(x)\equiv 1$ (see Remark \ref{remark: general slowing chars}).
Then, we get an explicit expression for the characteristic variable $w_i$
\begin{equation}
    w_i(t,x) = \mathring{w}_i(x-\lambda_i t)  \exp\left(-|\lambda_i| \int_0^t d_i(x-\lambda_i(t-r)) \, dr\right),
\end{equation}
which, taking the change of variables $z=x-\lambda_i(t-r)$, becomes
\begin{equation}
        \label{eq: explicit solution  SDO characteristic no slowing}
    w_i(t,x) = \mathring{w}_i(x-\lambda_i t)  \exp\left(-\mathrm{sgn}(\lambda_i)\int_{x-\lambda_i t}^{x} d_i(z) \, dz\right).
\end{equation}
Consequently, as long as $x<x_s$, the information of each characteristic field
is transported without modification and 
it decays exponentially at a non-local rate inside the sponge layer.

In \cite{karni_far-field_1996}, it is proposed to leave the left-going waves unchanged, i.e., to set 
\begin{align}
    \label{eq: not modifying left-going slowing damping}
        d_i(x) = 0,\; s_i(x) = 1,  \quad &\text{if } \lambda_i < 0.
\end{align} 
This way, physical waves that are expected to propagate back into the computational
domain remain unperturbed. 
In the linear case, if no information is in the sponge layer at the initial time,
there is nothing to be propagated back into the computational domain
through the left-going waves, so the choice of $d_i$ and $s_i$ for $\lambda_i < 0$
does not affect the solution.
In the nonlinear case, however, different characteristic fields might 
interact with each other
in the sponge layer, and the choice \eqref{eq: not modifying left-going slowing damping}
could allow to preserve physical incoming information,
as will be demonstrated in Section {\ref{sec: performance comparison}}.
In the following, we will use the abbreviation SDO to refer to
the modified system \eqref{eq: linear SDO} with 
\eqref{eq: definition d and s SDO}, and  \eqref{eq: not modifying left-going slowing damping}.

Alternatively, one could choose to damp not only the outgoing waves, but all the characteristic fields
inside the sponge layer.
If given $d(x)$ and $s(x)$ defined as in \eqref{eq: definition d and s SDO},
we take $d_i(x)=\frac{d(x)}{|\lambda_i|}$ and $s_i(x)=s(x)$
for \emph{every} $i\in\{1,\ldots,N\}$,
from \eqref{eq: linear SDO} we recover the familiar system
\begin{subequations}
    \label{eq: scalar SDO}
    \begin{align}
    \partial_t\bfq + s(x) 
    \bfA \partial_x\bfq  
    &= - \bfR \diag\left(d(x) s(x)\right) \bfR^{-1} \left(\bfq-\overline{\bfq}\right)\\
    &= -d(x) s(x) 
    \left(\bfq-\overline{\bfq}\right).
\end{align}
\end{subequations}
Damping source terms with a similar form
have been applied previously to, e.g., the wave equation \cite{israeli_approximation_1981}, the 
linearized
Euler equations \cite{goodrich_comparison_2008}, 
and the Schr\"odinger equation \cite{kosloff_absorbing_1986}.
For the sake of comparison, we will also consider 
\eqref{eq: scalar SDO}
and refer
to it as scalar slowing-down and damping operators (S-SDO).

Even though the SDO technique was introduced for a linear system,
a natural way to use it for nonlinear systems is
to compute the far-field terms \eqref{eq: definition far-field operators} (now dependent on $\bfq$)
locally \cite{karni_accelerated_1992}.
This will be the approach considered
in Section \ref{sec:Results} for the
nonlinear equations \eqref{eq:conservative form lagrangian}.
Additionally, an alternative extension to the nonlinear setting is proposed in Section \ref{sec:NDO}.

\begin{remark}
\label{remark: general slowing chars}
In the general case where $s_i$ is not constant, at any given time $t$ and position $x$, 
it is possible to track the characterstic curves
$x=X(t)$ back to $x_0$. It suffices to solve \eqref{eq:implicit characteristic}, which 
yields $t=\int_{x_0}^x \frac{1}{\lambda_i s_i(z)} dz$.
Then, by defining $T(x)=\int_{-\infty}^x \frac{1}{\lambda_i s_i(z)} dz$
and noting that $T$ is invertible in $[0,x_\infty[$,
we can write $x=T^{-1}(t+T(x_0))$ and $x_0=T^{-1}(T(x)-t)$.
\end{remark}

\subsection{The relaxation method (RM)}
\label{sec:RM}

The relaxation method
 is a postprocessing technique where the numerical solution is forced 
 towards a prescribed function
after each time step.
The process
is described in Algorithm \ref{alg: relaxation scalar}, where the so-called weight function $\Gamma(x)$ decreases
from $1$ to $0$ in the sponge layer, as depicted in 
Figure \ref{fig: Shape slowing damping funcs}.
In the relaxation step, 
$\bfq$ is replaced by
a convex combination of the numerical solution $\bfq^*$ and some target function, which in this case is
the constant far-field state $\overline{\bfq}$. 
    \begin{algorithm}
    \begin{algorithmic} \caption{Relaxation method with scalar weight function} \label{alg: relaxation scalar}
    \Require $\qn$
    \Ensure $\qn[n+1]$
    \State $\bfq^* \gets $ Evolve 
    a discretization of \eqref{eq:conservation law} over a time step $\Delta t$ with
    \\ \quad initial data $\qn$ \Comment{Hyperbolic step}
    \State $\qn[n+1] \gets \Gamma(x)\bfq^*+(1-\Gamma(x))\overline{\bfq}$ \Comment{Relaxation step}
    \end{algorithmic}
\end{algorithm}

It is well known that the performance of this ABC depends on the choice of the weight function \cite{choi_performance_2020}. 
By considering the method as a nodal filtering
procedure, Engsig-Karup \cite{engsig-karup_unstructured_nodate}
shows that
a necessary condition to preserve the continuity of the modal coefficients through
elements is that
\begin{equation}
    \Gamma(x_s)=1, \quad \Gamma'(x_s)=\Gamma''(x_s)=\Gamma'''(x_s)=...=0,
\end{equation}
which is similar to the smooth transition requirement on the slowing-down and damping functions proposed by Karni \cite{karni_far-field_1996}.
Based on this observation, in \cite{engsig-karup_nodal_2006}
the weight function 
\begin{equation}
    \label{eq: Engsig-Karup weight function}
    \Gamma_A(x)=-2\left(1-\phi(x)\right)^3+3\left(1-\phi(x)\right)^2,
\end{equation}
is proposed, where
\begin{equation}
    \label{eq:scaling weight functions}
    \phi(x)=\left(\frac{x-x_s}{x_\infty-x_s} \right)\chi_{]x_s,x_\infty]}(x)+\chi_{]x_\infty,+\infty[}(x).
\end{equation}
In \cite{mayer_fractional_1998}, Mayer et al. propose the weight function
\begin{equation}
    \label{eq: Mayer weight function}
    \Gamma_B(x) = 1-\left[b\phi(x)^3+(1-b)\phi(x)^6\right],
\end{equation}
where $b\in [0,1]$ is a parameter that determines its rate of decrease.

The accuracy of this method also depends on
the sponge layer length.
A customary rule of thumb to design sponge layers
for the relaxation method is the 
requirement $|x_\infty-x_s|>2L$, where $L$ is the largest wavelength
present in the solution \cite{mayer_fractional_1998, engsig-karup_nodal_2006, choi_performance_2020}.
While this restriction might be inconvenient for problems involving 
long waves, it is not unique to the relaxation method.
For instance, Perić and Abdel-Maksoud \cite{peric_reliable_2016} propose the same requirement on the 
sponge layer for their linear and quadratic damping  far-field operators.

In order to improve the properties of this method, a popular approach is to
define a specialized weight function.
Hu et al. \cite{hu_numerical_2015} choose the weight function at each
time step such that the mass absorbed artificially is equal to the
mass entering the sponge layer. 
They report that, with this approach, acceptable results are attained with the looser restriction
$|x_\infty-x_s|>0.4 L$
for the 2D Navier-Stokes equations
with free surface.
In \cite{seng_slamming_2012}, Seng notes that the damping
effect of the relaxation method increases dramatically as the time step
size is refined.
To mitigate this phenomenon, the author proposes a weight function that is
updated
dynamically
depending on the distance between the numerical and target solutions.

   A variation of the relaxation method was proposed by Wasistho et al. \cite{wasistho_simulation_1997}
for the 2D Navier-Stokes equations, where the authors proposed the application of relaxation
after each stage of the time integration with a Runge-Kutta scheme. 
Some numerical results concerning this approach are presented in
Section \ref{sec: performance comparison}.

\section{Extensions of some sponge layer-based ABCs}
\label{sec: Extensions}

\subsection{A relation between the RM and the far-field damping operator}
\label{sec:relaxation and far-field}
By performing some
manipulations to the relaxation step
in Algorithm \ref{alg: relaxation scalar},
we
obtain
\begin{subequations}
\begin{align}
    \qn[n+1]&=\Gamma(x)\bfq^*+(1-\Gamma(x))\overline{\bfq} \label{eq: Relaxation step}\\
    &=\bfq^*+\Delta t \frac{1-\Gamma(x)}{\Delta t}(\overline{\bfq}-{\bfq}^*)\\
    &=: \bfq^* +\Delta t \,\Tilde{\Gamma}(x;\Delta t)(\overline{\bfq}-{\bfq}^*), 
\end{align}
\end{subequations}
which is a forward Euler step to approximate the solution
 $\bfq(\Delta t)$ to the system
\begin{equation}
    \label{eq: relaxation forward euler scalar}
    \partial_t \bfq = \Tilde{\Gamma} (x;\Delta t)(\overline{\bfq}-\bfq),\quad \bfq(0)=\bfq^*.
\end{equation}
Hence, Algorithm \ref{alg: relaxation scalar} can be seen as
a discretization 
of
\begin{equation}
    \label{eq: relaxation far-field 1st order}
    \partial_t \bfq + \partial_x \bff(\bfq) = \Tilde{\Gamma} (x;\Delta t)(\overline{\bfq}-\bfq),
\end{equation}
using Lie-Trotter operator splitting.

The right hand side of Equation \eqref{eq: relaxation far-field 1st order} can be written as
\begin{subequations}
\begin{align}
    \Tilde{\Gamma}(x;\Delta t)(\overline{\bfq}-\bfq) 
    &= 
    \Tilde{\Gamma}(x;\Delta t)  \bfR \bfR^{-1}(\overline{\bfq}-\bfq) \\
    &=\bfR \begin{pmatrix}
 \frac{\Tilde{\Gamma}(x;\Delta t) |\lambda_1|}{|\lambda_1|}  & 0 & \cdots & 0 \\
0 &\frac{\Tilde{\Gamma}(x;\Delta t) |\lambda_2|}{|\lambda_2|}  & \ddots & \vdots \\
\vdots & \ddots & \ddots & 0 \\
0 & \cdots & 0 & \frac{\Tilde{\Gamma}(x;\Delta t) |\lambda_N|}{|\lambda_N|} 
\end{pmatrix}
\bfR^{-1} (\overline{\bfq}-\bfq) \\
&= \bfR {\boldsymbol{\Lambda}_{\bf D}}(x;\Delta t)  \bfR^{-1}(\overline{\bfq}-\bfq)=\bfD(x;\Delta t) (\overline{\bfq}-\bfq),
\end{align}
\end{subequations}
using the notation from Section \ref{sec: SDO}.
It follows that, with the RM, the $i$-th characteristic field
is absorbed inside the sponge layer at the
damping rate
\begin{equation}
    \label{eq: damping rate scalar relaxation}
    d_i^{\text{RM}}(x;\Delta t) = \frac{\Tilde{\Gamma}(x;\Delta t)}{|\lambda_i|}=\frac{1-\Gamma(x)}{\Delta t |\lambda_i|}.
\end{equation}
The inversely proportional relation between $d_i$ and $\Delta t$ 
is consistent with the increase in damping severity as the time step size decreases using the RM
observed
in the literature (see, e.g., \cite{seng_slamming_2012}).

\subsection{The RM with matrix-valued weight function (RM-M)}
\label{sec: relaxation matrix valued}

We propose a reframing of the RM
such that it resembles the directional behavior of Karni's
damping far-field operator.
We achieve this by allowing the weight function $\bGamma$ to be matrix-valued and defined 
in terms of $\bfR$ (and therefore $\bfq$ in the nonlinear case).
The process is described in Algorithm \ref{alg: relaxation matrix} and justified as follows.
Proceeding analogously to the previous section, the relaxation step in Algorithm \ref{alg: relaxation matrix}
\begin{equation}
\label{eq: Relaxation step Matrix}
    \qn[n+1] = {\bGamma}(x)\bfq^*+(\bfI-\bGamma(x))\overline{\bfq},
\end{equation}
can be seen as a forward Euler discretization of
\begin{equation}
    \label{eq: relaxation forward euler matrix}
    \partial_t \bfq = \Tilde{\bGamma} (x,\bfq;\Delta t)(\overline{\bfq}-\bfq),
\end{equation}
where
\begin{equation}
    \Tilde{\bGamma} (x,\bfq;\Delta t)=
    \frac{1}{\Delta t}(\bfI-{\bGamma}(x))=
     \bfR 
     \begin{pmatrix}
 \frac{1-\Gamma_1(x)}{\Delta t} & 0 & \cdots & 0 \\
0 &\frac{1-\Gamma_2(x)}{\Delta t} & \ddots & \vdots \\
\vdots & \ddots & \ddots & 0 \\
0 & \cdots & 0 & \frac{1-\Gamma_N(x)}{\Delta t}
\end{pmatrix}
    \bfR^{-1}.
\end{equation}

Thereby, up to first order of accuracy in time,
we approximate a system of the form \eqref{eq: linear SDO},
where the desired outgoing wave is damped at the rate
\begin{equation}
    \label{eq: damping rate matrix relaxation}
    d_i^{\text{RM-M}}(x) = \frac{1-\Gamma_i(x)}{\Delta t |\lambda_i|},
\end{equation}
while the other waves remain unperturbed.
The relaxation method of Section 
\ref{sec:RM} is a particular case of this formulation, where \emph{all} 
the waves are modified at the non-zero rate \eqref{eq: damping rate matrix relaxation} inside the sponge layer.

If $\bfR$ is independent of $\bfq$, the matrix-valued weight function ${\boldsymbol{\Gamma}}$ 
does not have to be computed at each time step in Algorithm \ref{alg: relaxation matrix}.
Additionally, if the far-field state $\overline{\bfq}$ is constant, we
can define the translated function $\bfv=\bfq-\overline{\bfq}$ that satisfies the conservation law
$\partial_t \bfv+\partial_x \bff(\bfv+\overline{\bfq})=0$.
This system may be solved with the far-field state $\bf0$
and $\bfq$ recovered afterward.
In this case, we have
\begin{equation}
    ||\bfv^{[n+1]}||_{2} \leq ||{\boldsymbol{\Gamma}}\bfv^*||_2 \leq ||{\boldsymbol{\Gamma}}||_2||\bfv^*||_2\leq ||\bfv^*||_2,
\end{equation}
thus the stability of the relaxation step is guaranteed.

\begin{algorithm}
    \begin{algorithmic} \caption{Directional relaxation method with matrix-valued weight function } \label{alg: relaxation matrix}
    \Require $\qn$
    \Ensure $\qn[n+1]$
    \State $\bfq^* \gets $ Evolve 
    a discretization of \eqref{eq:conservation law} over
    a time step $\Delta t$ with\\
    \quad initial data $\qn$ \Comment{Hyperbolic step} 
    \For{each characteristic field $(\lambda_i, \bfr_i)$}
    \If{$i$-th wave is being damped}
    \State $\Gamma_i(x) \gets \Gamma(x)$
    \Else
    \State $\Gamma_i (x)\gets 1$
    \EndIf
    \EndFor
    \State ${\boldsymbol \Lambda}_{\boldsymbol{\Gamma}}(x) \gets $diagonal matrix-valued  function with entries $\Gamma_i(x)$
    \State ${\boldsymbol {\Gamma}}(x) \gets \bfR {\boldsymbol \Lambda}_{\boldsymbol{\Gamma}}(x) \bfR^{-1}$  where $\bfR$ contains the right eigenvectors \\
    \quad  of $\bff ' (\bfq)$
    \State $\qn[n+1] \gets {\boldsymbol\Gamma}(x)\bfq^*+(\bfI-{\boldsymbol \Gamma}(x))\overline{\bfq}$ \Comment{Relaxation step}
    \end{algorithmic}
\end{algorithm}

\begin{figure}
    \centering
    \includegraphics[width=\textwidth]{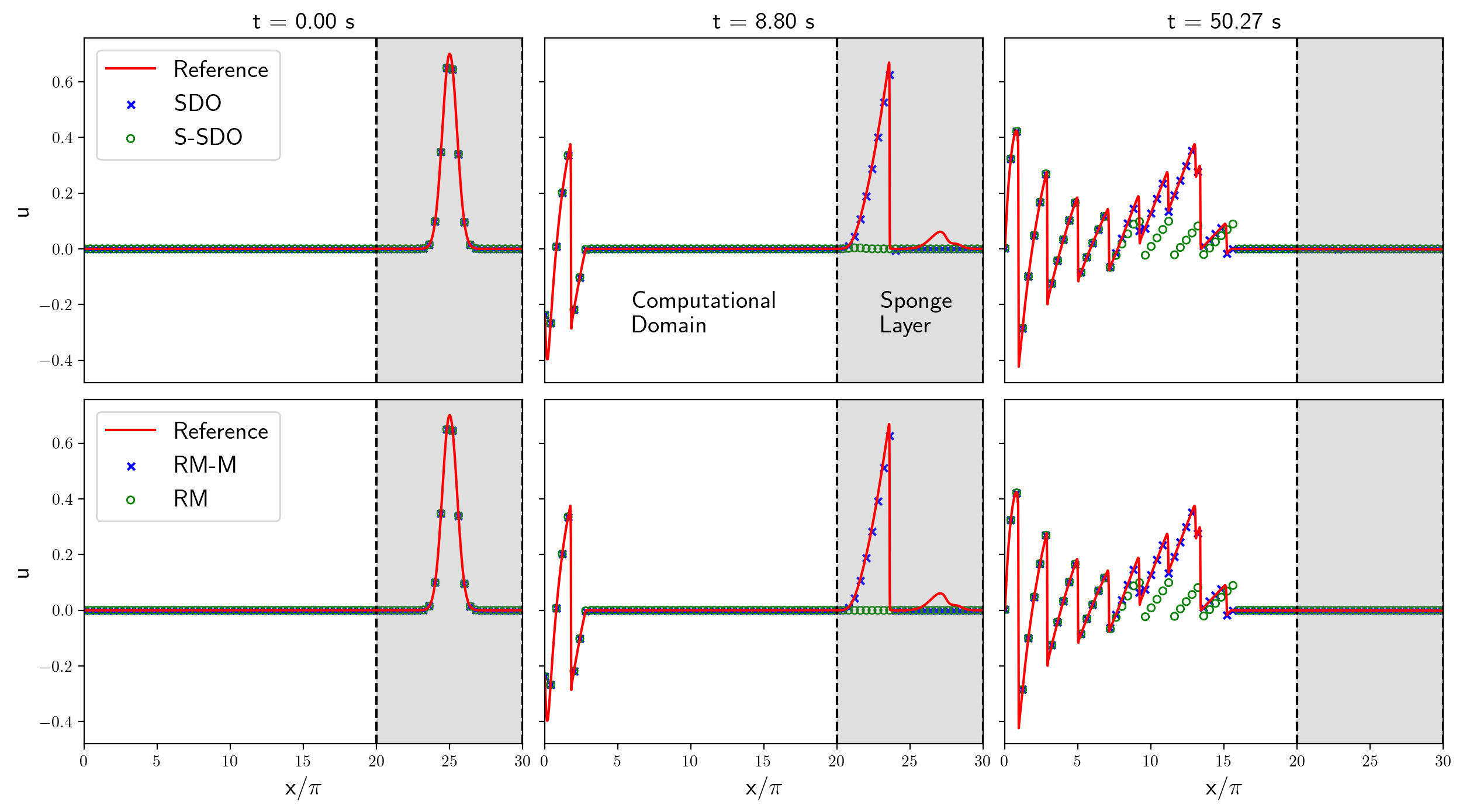}
    \caption{
    A Gaussian pulse is initialized inside the sponge layer to represent the effect of nonlinear interactions 
    taking place there. Thus, we consider its left-going component as a physical wave to be preserved.
    The right-going sawtooth-like wave corresponds to disturbances introduced by the motion of a piston, which will
    be described in Section \ref{sec:piston problem}.
    The sponge layer is colored light grey.
    Above: Effect of directional (SDO) and scalar (S-SDO) damping operators applied to the incoming physical wave.
    Below: Effect of directional relaxation method
    with matrix-valued weight function  (RM-M)
    and scalar relaxation method (RM) applied to the incoming physical wave.}
    \label{fig: artificial problem}
\end{figure}

We illustrate the effect of damping waves corresponding to a particular characteristic field
with a preliminary numerical experiment using the framework to be described in Section \ref{sec:Numerics}.
An artificial setup where an incoming pulse is created inside the sponge layer is considered,
as shown in Figure \ref{fig: artificial problem}.
The initial pulse is immediately damped by the
scalar methods (RM and S-SDO), while just the right-going
components of it are absorbed by the directional methods (RM-M and SDO).
In a practical setting, the incoming pulse could correspond to, e.g.,
physical nonlinear interactions inside the sponge layer.
While an artificially large pulse is used in Figure \ref{fig: artificial problem} for visualization purposes,
such nonlinear interactions might not be less significant than spurious reflections due to the boundary treatment (see Section \ref{sec: performance comparison}).

\subsection{On attaining higher order of accuracy in time with the RM}
\label{sec: Strang}

Let us consider Algorithm \ref{alg: relaxation matrix} for the linear problem \eqref{eq:conservation law linear},
where the matrix-valued weight function $\boldsymbol{\Gamma}$ is independent of $\bfq $.
Then, the relaxation step \eqref{eq: Relaxation step Matrix}
could be written as a high-order discretization of a system 
$\partial_t \bfq = \mathcal{R}(g(\boldsymbol{\Gamma}(x),\Delta t))\, (\overline{\bfq}-\bfq)$, where
$\mathcal{R}$ is the stability function of some
Runge-Kutta scheme \cite{gottlieb_high_2009}.
Moreover, we can interpret it as the \emph{exact} solution to the problem
\begin{equation}
\label{eq: target problem higher-order splitting}
    \partial_t \bfq = \bfW(x;\Delta t)(\overline{\bfq}-\bfq),\quad \bfq(0)=\bfq^*,
\end{equation}
where $\bfW$ satisfies
\begin{equation}
    {\boldsymbol\Gamma} (x;\Delta t)= \exp({-\Delta t \bfW(x;\Delta t)}).
\end{equation}
Since ${\boldsymbol\Gamma}= 
\bfR {\boldsymbol \Lambda}_{\bGamma} \bfR^{-1}$, then $\bfW=\bfR {\bm{\Lambda}_{\bf W}} \bfR^{-1}$, where
\begin{equation}
    \label{eq: damping with singularity 1}
      w_i(x) = 
      \frac{-\ln\left(\Gamma_i(x)\right)}{\Delta t},
\end{equation}
with $\Gamma_i$ and $w_i$ denoting the main diagonal entries of the matrices
${\boldsymbol \Lambda}_{\boldsymbol{\Gamma}}$ and ${\boldsymbol \Lambda_W}$ respectively.
Thus, the $i$-th wave is damped at the rate
\begin{equation}
    \label{eq: log damping rate}
    d_i^{\bfW}(x) = \frac{w_i(x)}{|\lambda_i|} =  \frac{-\ln\left(\Gamma_i(x)\right)}{\Delta t |\lambda_i|},
\end{equation}
and will be completely absorbed as $\Gamma_i\to 0$,
i.e., at the boundary of the sponge layer $x_\infty$.
This further justifies
the approaches proposed by Seng \cite{seng_slamming_2012} and
Wasistho et al. \cite{wasistho_simulation_1997}
where, through different arguments, $\Gamma$ is substituted by 
$\Gamma^{C \Delta t}$ to reduce the damping sensitivity to the time step size.

As a consequence of the Lie-Trotter splitting, a first-order error in time can be expected when
approximating a solution to 
the system
\begin{equation}
    \label{eq: relaxation far-field 2nd order}
    \partial_t \bfq + \partial_x \bff(\bfq) = \bfW(x;\Delta t)(\overline{\bfq}-\bfq),
\end{equation}
and a natural approach to get a higher-order accurate in time approximate solution to \eqref{eq: relaxation far-field 2nd order}
would be to use the RM in the
context of high-order operator splitting.
For instance, we could use Strang splitting as described in Algorithm \ref{alg: Strang splitting}.
However, in this setting, both splitting techniques are equivalent up to the choice of the scalar weight functions $\Gamma_i$.
This can be seen as follows:
Assuming that the solution is initialized as $\overline{\bfq}$ inside the sponge layer,
the
first relaxation step at $t^{[0]}=0$ does not alter the solution.
Similarly, the second
relaxation step in Algorithm \ref{alg: Strang splitting}
at the final time
does not affect the solution in the computational domain,
since it is local and no hyperbolic step is taken afterward.
Finally, let $\bfq^{**}(t^{[n]})$ denote the second stage of Algorithm \ref{alg: Strang splitting} with
initial data $\bfq^{[n]}$, and let $\bfq^{*}(t^{[n+1]})$ denote the first stage of Algorithm \ref{alg: Strang splitting} with
initial data $\bfq^{[n+1]}$.
Then, the following relation holds:
\begin{subequations}
\begin{align}
\label{eq: equivalence Lie-Trotter Strang linear}
    \bfq^{*}(t^{[n+1]}) &= \boldsymbol{\Gamma} \bfq^{[n+1]} +(\bfI - \boldsymbol{\Gamma})\overline{\bfq} \\
    &=\boldsymbol{\Gamma} \left(\boldsymbol{\Gamma} \bfq^{**}(t^{[n]})+(\bfI - \boldsymbol{\Gamma})\overline{\bfq}\right)+(\bfI - \boldsymbol{\Gamma})\overline{\bfq} \\
    &= \boldsymbol{\Gamma}^2 \bfq^{**}(t^{[n]})+(\bfI - \boldsymbol{\Gamma}^2)\overline{\bfq}.
\end{align}
\end{subequations}
The latter amounts to taking a single relaxation step,  i.e., using Algorithm \ref{alg: relaxation matrix},
with the 
modified scalar weight functions $\Gamma_i^2$.
Therefore,  in the linear case,
as long as the scalar weight functions $\Gamma_i$ are tuned
appropriately 
for the application at hand, 
a second-order in-time approximation of a system of the form
\eqref{eq: relaxation far-field 2nd order} is obtained with
Algorithm \ref{alg: relaxation matrix}.

On the other hand, since $\boldsymbol{\Gamma}$ depends on $\bfq$
in the nonlinear case, the last observation is no longer valid, and
Algorithms \ref{alg: relaxation matrix} and \ref{alg: Strang splitting} are not equivalent a priori. 
The effect of applying the RM as in Algorithm \ref{alg: Strang splitting}
in the nonlinear setting is evaluated
numerically in Section \ref{sec: performance comparison}.

\begin{algorithm}
    \begin{algorithmic} \caption{Relaxation method with Strang splitting } \label{alg: Strang splitting}
    \Require $\qn$
    \Ensure $\qn[n+1]$
    \State Compute the weight function ${\boldsymbol{\Gamma}}(x)$ as in Algorithm \ref{alg: relaxation matrix}
    \State $\bfq^{*} \gets {\boldsymbol\Gamma}(x)\qn+(\bfI-{\boldsymbol \Gamma}(x))\overline{\bfq}$  \\
    \quad $\left(\bfq (\Delta t/2) \text{ solving  \eqref{eq: target problem higher-order splitting} with }\bfq(0)=\qn\right)$
    \Comment{Relaxation step}
    \State $\bfq^{**}\gets $ Evolve 
    a discretization of \eqref{eq:conservation law} over
    a time step $\Delta t$ with\\
    \quad initial data $\bfq^{*}$ \Comment{Hyperbolic step}
    \State $\qn[n+1] \gets {\boldsymbol\Gamma}(x)\bfq^{**}+(\bfI-{\boldsymbol \Gamma}(x))\overline{\bfq}$ \\
    \quad $\left(\bfq (\Delta t/2) \text{ solving  \eqref{eq: target problem higher-order splitting} with }\bfq(0)=\bfq^{**}\right)$
    \Comment{Relaxation step}
    \end{algorithmic}
\end{algorithm}

\subsection{Nonlinear damping far-field operator (NDO)}
\label{sec:NDO}
We propose an alternative extension of the damping operator 
from Section \ref{sec: SDO} to hyperbolic systems of the 
form \eqref{eq:conservation law} with a general nonlinear
flux function.
Here, the change of variables $\bfw = \bfR^{-1}\bfq$
no longer allows us
to decompose a nonlinear system of the form 
\eqref{eq:conservation law} into a system of scalar 
evolution equations.
 Instead,
we consider
the nonlinear characteristic variables $\boldsymbol{\alpha} = \bfR^{-1}\partial_x \bfq$ which, after some manipulations to System \eqref{eq:conservation law} (see \cite[Ch.1]{lefloch_hyperbolic_2002}),
yield
\begin{equation}
    \label{eq: decoupled nonlinear}
    \partial_t \alpha_i + \partial_x \left( \lambda_i(\bfq) \alpha_i\right) = 
    \mathcal{I}(\bfq,\boldsymbol{\alpha}),
    \quad 1\leq i \leq N,
\end{equation}
where $\mathcal{I}(\bfq,\boldsymbol{\alpha})$ denotes the nonlinear interactions between the $N$
characteristic fields.
Although this is not a decoupled system of equations, we can add a damping term
at the right-hand side, resembling the modification
done in \eqref{eq:each characteristic} for the linear equations \eqref{eq:conservation law linear}, to obtain
\begin{equation}
    \partial_t \alpha_i + \partial_x \left( \lambda_i(\bfq) \alpha_i\right) =
    \mathcal{I}(\bfq,\boldsymbol{\alpha})-d_i(x) |\lambda_i(\bfq)| \alpha_i,\quad 1\leq i \leq N.
\end{equation}
After multiplying this system from the left by $\bfR$ and integrating from $x$ to $+\infty$ (inverting the process carried out to obtain \eqref{eq: decoupled nonlinear}), we get
\begin{subequations}
\label{eq: NDO}
\begin{align}
    \partial_t \bfq(t,x) + \partial_x \bff(\bfq(t,x)) & = \int_x^{+\infty} \bfR
    \boldsymbol{\Lambda}_\bfD(z){\bfR}^{-1} \partial_z \bfq (t,z)\,dz\\
    & \approx  \int_x^{x_\infty} \bfD(z,\bfq) \partial_z \bfq (t,z)\,dz, \label{eq: NDO integral damping term}
\end{align}
\end{subequations}
where the artificial boundary at the end of the sponge layer is enforced through the
assumption  $\bfq(t,x)=\overline{\bfq}$ for $x\geq x_\infty$.
If the damping rate functions $d_i$ are defined as in Section
\ref{sec: SDO}, the term \eqref{eq: NDO integral damping term}
vanishes inside of the computational domain $x<x_s$.
Moreover, if the flux is linear and each $d_i$ is constant,
the damping far-field operator of Section \ref{sec: SDO} is
recovered. 

Some challenges that arise when dealing with
the nonlinear damping term \eqref{eq: NDO integral damping term}
are its non-local nature and the presence of a
nonconservative product.
This term is significantly simplified
if one can write
$\bfD(z,\bfq)=\partial_\bfq \bfh(z,\bfq)$
for some function
$\bfh\in C^1(\real\times\real^N,\real^N)$.
Indeed, in that case, we would have
\begin{subequations}
\begin{align}
    \int_x^{x_\infty} \bfD(z,\bfq) \partial_z \bfq (t,z)\,dz &= 
    \int_x^{x_\infty} \partial_\bfq \bfh(z,\bfq)\partial_z\bfq(t,z)\,dz\\
    &=
    \int_x^{x_\infty} \frac{d}{dz} \bfh(z,\bfq) -  \partial_z\bfh(z,\bfq)\,dz\\
    &= \bfh(x_\infty,\overline{\bfq})-\bfh(x,{\bfq(t,x)})-\int_x^{x_\infty} \partial_z\bfh(z,\bfq(t,z))\,dz. \label{eq: NDO simplified}
\end{align}
\end{subequations}
Nevertheless, 
the integrability conditions to guarantee 
the existence of such a function $\bfh$
do not hold for a general system and choice of $d_i$'s
\cite{papachristou_aspects_2019}.
Thus, for the application studied in Section \ref{sec:Results} we work with a  discretization of the term \eqref{eq: NDO integral damping term}, described in Appendix \ref{apen: Discretization NDO}.

\section{The Euler equations}
\label{sec:Euler equations}

As an application for the techniques presented in Sections \ref{sec:Boundary treatment} and \ref{sec: Extensions}, we consider the one-dimensional formulation
of the Euler equations for an inviscid compressible gas 
\begin{subequations}
\label{eq:Euler 1D}
    \begin{align}
        \partial_t \rho + \partial_{\xi} \left(\rho u\right) &= 0, \label{eq:cons mass 1D}\\
        \partial_t \left(\rho u\right) +  \partial_{\xi} \left( \rho u^2+p\right) &=0 ,\label{eq:cons momentum 1D}\\
        \partial_t \barE+ \partial_{\xi} \left[ (\barE+p) u \right] &= 0.\label{eq:cons energy 1D}
    \end{align}
\end{subequations}

Here $\rho$ denotes the mass density, $u$ the flow velocity, $p$ the pressure, and $\rho u$ the momentum.
The energy density given by $ \barE=\frac12 \rho u^2 +\rho e$,
where $\rho u^2$ is the kinetic energy and $e$ is the specific internal energy.
Equations \eqref{eq:cons mass 1D}, \eqref{eq:cons momentum 1D}, and \eqref{eq:cons energy 1D} describe the 
conservation of mass, momentum, and energy respectively, and
can be derived from the Navier-Stokes equations or foundational principles in
continuum mechanics if viscosity and heat flux are neglected \cite{landau_ideal_1987}. 

We deal with an ideal gas with constant specific heats, hence we can define
the constant heat capacity ratio $\gamma=c_p/c_V$, and
our equation of state becomes \cite{whitham_linear_1974}
\begin{equation}
    \label{eq:EOS}
    e = \frac{p}{(\gamma-1)\rho}.
\end{equation}
Throughout this work, we assume the density $\rho$ to be strictly positive for all times.
Therefore vacuum states, which require special analytical 
and numerical treatment \cite[Sec. 4.6]{toro_riemann_2009},
are not considered.

It is often convenient
to write Equations \eqref{eq:Euler 1D} in Lagrangian form.
Given a real function $\tilde {\xi}(\cdot)$, we will keep track of a set of conserved quantities at each particle defining the trajectory map \cite{phan_thi_my_finite_2021}
    \begin{equation}
    \label{eq:trajectory map}
        x(\xi,t) = \int_{\tilde{\xi}(t)}^\xi \rho(\tau,t) d\tau.
    \end{equation}
From our non-vacuum state assumption $\rho$ is strictly positive,
hence $x(\cdot\,,t)$ is bijective and we can use the correspondence $f(x(\xi,t),t)=f(\xi,t)$.
After defining the Lagrangian specific volume 
\begin{equation}
    V(x(\xi,t),t)=\frac{1}{\rho(\xi,t)},
\end{equation}
and the material time derivative operator (with a slight abuse of notation)
\begin{equation}
    \partial_t := \partial_t + u\,\partial_\xi,
\end{equation}
the equations \eqref{eq:Euler 1D} become 
\begin{subequations}
\label{eq:Lagrangian Euler 1D}
    \begin{align}
        \partial_t V - \partial_x u = 0,\\
        \partial_t u + \partial_x p = 0,\\
        \partial_t E + \partial_x (up) = 0,
    \end{align}
\end{subequations}
where $E(x,t)=\barE/\rho$ \cite{smoller_shock_1994}.
Furthermore, using the equation of state \eqref{eq:EOS}, we have the 
closure relation
\begin{equation}
    \label{eq:relation p and V}
    p=(\gamma-1)\left(E-\frac12 u^2\right)/V.
\end{equation}

The hyperbolic system \eqref{eq:Lagrangian Euler 1D} can be more compactly written as
\begin{equation}
    \label{eq:conservative form lagrangian}
    \partial_t \bfq + \partial_x \bff(\bfq) = {\bf 0},
\end{equation}
where 
$
    \bfq = \left(q_1,q_2,q_3\right)^T
    =\left(V,u,E\right)^T
$
and
$
    \bff(\bfq)= \left(q_1,p(\bfq),q_2 p(\bfq) \right)^T
    = \left(-u,p,up\right)^T
$

 denote the vector of conserved quantities and the nonlinear flux function respectively.
For smooth solutions, System \eqref{eq:conservative form lagrangian} can be written in quasilinear form
\begin{equation}
    \label{eq:quasilinear form lagrangian }
    \partial_t\bfq+\bff '(\bfq) \partial_x \bfq = {\bf 0},
\end{equation}
where the Jacobian matrix
\begin{equation}
{\bf f}'({\bf q}) =
\left(\begin{array}{ccc}
0 & -1 & 0 \\
-\frac{p}{V} & \frac{{u}(1-\gamma)}{{V}} & \frac{-1+\gamma}{{V}} \\
-\frac{{p} {u}}{{V}} & p+\frac{{u}^2(1-\gamma)}{{V}} & \frac{{u}(-1+\gamma)}{{V}}
\end{array}\right),
\end{equation}
has symmetric eigenvalues 
$
\lambda_1(\bfq)=-\sqrt{\gamma p/V}$, $\lambda_2=0$, $ \lambda_3(\bfq)=\sqrt{\gamma p/V}$,
with corresponding eigenvectors
$
\boldsymbol{r}_{1,3}(\boldsymbol{q}) = \left( -\frac{\sqrt{V}}{\sqrt{p}}, \mp\sqrt{\gamma}, \sqrt{pV}\mp u\sqrt{\gamma} \right)^T,
\,
\boldsymbol{r}_2(\boldsymbol{q}) = \left( \gamma-1, 0, p \right)^T.
$

It can be shown that the 
characteristic fields corresponding to
the eigenvalues $\lambda_1$ and $\lambda_3$ are genuinely nonlinear while the one
corresponding to $\lambda_2$ is linearly degenerate.
Furthermore, from our non-vacuum state assumption, it follows that $p>0$ at all times, 
thus the system is strictly hyperbolic.

We also consider a linearization  of System \eqref{eq:conservative form lagrangian} around the state $\overline{\bfq}=(1, 0, \gamma/(\gamma-1))$ 
\begin{equation}
\label{eq:linearized equations}
    \partial_t\bfq+\bfA \partial_x \bfq={\bf 0},
\end{equation}
where the matrix
\begin{equation}
    \bfA = \begin{pmatrix}
        0 & -1 & 0\\
        -\gamma & 0 & \gamma-1\\
        0 & \gamma & 0
    \end{pmatrix},    
\end{equation}
has eigenvalues $\lambda_1=-\gamma$, $\lambda_2=0$, $\lambda_3=\gamma$,
and eigenvectors
$
\boldsymbol{r}_1 = \left(-\frac{1}{\gamma}, -1, 1 \right)^T,$
$ \boldsymbol{r}_2 = \left( -\frac{1-\gamma}{\gamma},  0,  1 \right)^T,$ 
$\boldsymbol{r}_3 = \left( -\frac{1}{\gamma},  1,  1 \right)^T.
$
 
\subsection{The piston problem}
\label{sec:piston problem}

Let us formulate the piston problem of gas dynamics \cite{inoue_propagation_1993}, which will be the main application 
for the numerical treatment of boundary conditions studied in this work.

The equations \eqref{eq:Euler 1D} are solved in the domain $\Omega_t=[\tilde{\xi}(t),+\infty[$
with a constant steady state 
    $\overline{\bfq} = (1, 0, \gamma/(\gamma-1))$
as initial conditions.
The left physical boundary $\tilde{\xi}(t)$ corresponds to the position of an oscillatory plate given by
\begin{equation}
    \label{eq: piston position}
    \tilde{\xi}(t)=M \left(\cos(t)-1\right)\ \mathbbm{1}_{[\real^+]}(t),
\end{equation}
where $M$ is a constant that will determine the amplitude of the oscillations.
Since the particle velocity of the polytropic gas must coincide with the velocity of the plate at $\tilde{\xi}(t)$,
we have the inflow boundary condition
\begin{equation}
\label{eq:piston velocity}
    u(\tilde{\xi}(t),t)= {\tilde{\xi}}'(t) = -M \sin(t).
\end{equation}
Moreover, using the Lagrangian coordinates \eqref{eq:trajectory map}, we can solve the system of equations 
\eqref{eq:conservative form lagrangian} over the
time-independent domain $\Omega=[0,+\infty[$ and circumvent the difficulties arising from a moving boundary.

In the recent work \cite[Ch. 3]{phan_thi_my_finite_2021}, Phan Thi My investigates outflow boundary conditions
for this problem by replacing the Euler equations with a Burgers-type model in the sponge layer. 
The simple wave behavior reduces spurious reflections even when the sponge layer contains only a 
small number of cells (fewer than ten). This approach is attractive because it achieves efficient absorption with 
minimal computational overhead. However, its performance does not improve when the sponge layer is thickened, 
which limits its flexibility.

\section{Numerical results}
\label{sec:Results}

We analyze the performance of the ABCs introduced in the previous sections 
for the nonlinear equations \eqref{eq:conservative form lagrangian} and their linearized
version \eqref{eq:linearized equations}.
For this, we consider the setup described at the beginning of Section \ref{sec:Boundary treatment}.
Throughout this section, "performance" refers to the effectiveness of an ABC in minimizing
wave reflections, quantified by the reflection error \eqref{eq: Error ABC}.

\subsection{Discretization details}
\label{sec:Numerics}

For the numerical results presented in this work, we use 
a high-order wave propagation formulation of Godunov-type methods, available in the
Clawpack software package \cite{mandli_clawpack_2016}
and accessible through the Pyclaw framework \cite{ketcheson_pyclaw_2012}.
For a detailed description of these methods, the reader is referred to
the works \cite{ketcheson_high-order_2013,ketcheson_wenoclaw_2008}.

Denoting by $Q_i$ the vector of conserved quantities $\bfq$ averaged over 
the volume $[x_{i-1/2},x_{i+1/2}]$, our 
semidiscretization of Systems \eqref{eq:conservative form lagrangian} and \eqref{eq:linearized equations} 
reads
\begin{equation}
\label{eq:semidiscretization-general-PDE}
    \partial_t Q_i=-\frac{1}{\Delta x}\left(\mathcal{A}^{+} \Delta q_{i-\frac{1}{2}}+\mathcal{A}^{-} \Delta q_{i+\frac{1}{2}}+\mathcal{A} \Delta q_i\right).
\end{equation}
The fluctuations $\mathcal{A}^{\pm} \Delta q_{i-\frac{1}{2}}$ 
represent the net contribution of the left and right-going waves 
emanating 
from the solution of a Riemann problem at the interface $x_{i-1/2}$.
The initial data for each Riemann problem is given by a piecewise polynomial 
reconstruction of the solution at both sides of the interface.
To avoid spurious oscillations around strong shocks without using an extremely refined grid,
a piecewise linear reconstruction with the TVD minmod limiter is used.
An advantage of the semi-discretization \eqref{eq:semidiscretization-general-PDE} is that
it can be applied
to systems in nonconservative form, like \eqref{eq: linear SDO}, in a natural way \cite{bale_wave_2003,ketcheson_high-order_2013}.
Since the semidiscretization and numerical inflow boundary condition 
(see Equations \eqref{eq: BC piston nonlinear} and \eqref{eq: BC piston linear} below) formally attain second order of accuracy,
we use Heun's two-stage second-order accurate Runge-Kutta scheme to integrate System \eqref{eq:semidiscretization-general-PDE} in time.
A Courant number $C=0.8$ is used to compute $\Delta t$ at each time step 
under the standard CFL condition
$
    \Delta t = {C \Delta x}/{\max_{i,p} {|\lambda_p(Q_i)|}}.
$

To compute an approximate Riemann solution at an interface with left and right constant initial data $q_{L},\, q_{R}$,
we use an HLL Riemann solver \cite{harten_upstream_1983}.
The approximate speeds are defined to bound
the wave speeds from the exact Riemann solution from below and above,
\begin{equation}
 s^1=\min _p\left(\min \left(\lambda_p^L, \lambda^{R}_p\right)\right) = \min \left(\lambda_1^L, \lambda_1^R\right), \quad
 s^2=\max _p\left(\max \left(\lambda_p^L, \lambda^{R}_p\right)\right) = \max \left(\lambda_3^L, \lambda^{R}_3\right).
\end{equation}
Since the eigenvalues of our Jacobian matrix are symmetric, we have $s^1=-s^2$. 
Therefore, in this case, the HLL solver coincides with the Local Lax-Friedrichs (or Rusanov) solver.
Since the solutions we are focused on do not develop vacuum states,
the use of the HLL solver helps to avoid their spurious appearance.

We now turn our attention to the piston problem described in Section \ref{sec:piston problem},
which will be used as a benchmark to evaluate the performance of the different ABCs.
Let us denote the $i$-th cell of the discretized domain by $\mathcal{C}_i=[x_{i-1/2},x_{i+1/2}]$,
and the 
ghost cells at the left and right boundaries
by $\mathcal{C}_{-1}$ and $\mathcal{C}_{N_x}$ respectively.
We assume the piston (of infinitesimal width)
to be located at the interface between the cells $\mathcal{C}_{-1}$ and
$\mathcal{C}_0$ ($x=0$).
Using the 
flow conditions at the piston, formulated  in \cite{fazio_central_2010}, and their
discretization described in \cite{phan_thi_my_finite_2021}, we get the 
second-order accurate approximations
\begin{subequations}
\label{eq: BC piston nonlinear}
\begin{gather}
    u^{[n]}_{-1}=-2M\sin(t^{[n]})-u_0^{[n]}, 
    \quad p^{[n]}_{-1}=p_0^{[n]}-M\cos (t^{[n]})\Delta x,\\
    V^{[n]}_{-1}=V^{[n]}_{0}\frac{\gamma+\delta_p}{\gamma-\delta_p}, 
    \quad \delta_p=\frac{p^{[n]}_{0} - p^{[n]}_{-1}}{p^{[n]}_{0} + p^{[n]}_{-1}}, 
    \quad E_{-1}^{[n]}=\frac{1}{\gamma-1}p_{-1}^{[n]}V_{-1}^{[n]}+\frac12 (u_{-1}^{[n]})^2.
\end{gather}
\end{subequations}
Following an analogous derivation for the linearized equations \eqref{eq:linearized equations}, we obtain
\begin{subequations}
\label{eq: BC piston linear}
\begin{gather}
    u^{[n]}_{-1}= - 2M\sin(t^{[n]})-u^{[n]}_0,\quad p^{[n]}_{-1}=p^{[n]}_0 -M\cos(t^{[n]}) \Delta x, \\
    V^{[n]}_{-1}=\frac{1}{\gamma}(p^{[n]}_0-p^{[n]}_{-1}+\gamma V^{[n]}_0), \quad E^{[n]}_{-1}=\frac{1}{\gamma-1}(p^{[n]}_{-1}+\gamma V^{[n]}_{-1}),
\end{gather}
\end{subequations}
where we have the closure relation $p=\gamma V+(\gamma-1)E$ instead of \eqref{eq:relation p and V}.
At the right boundary we set $Q_{N_x} = \overline{\bfq}$, 
where $\overline{\bfq}$ is the far-field state
$\overline{\bfq
}=(1, 0, \gamma/(\gamma-1))$
for the nonlinear equations \eqref{eq:conservative form lagrangian},
and $\overline{\bfq}=\bf0$
for the linear equations \eqref{eq:linearized equations}. 
Here, $\mathcal{C}_{N_x}$ is the ghost cell with left boundary $x_{\infty}$.
For both systems, the initial condition will be
their corresponding far-field state $\overline{\bfq}$.

For all the numerical tests below, $x_s=20\pi$, and $x_\infty$
takes different values to evaluate the dependence of the ABCs on the sponge layer's length $\omega = |x_\infty-x_s|$.
Since the piston's frequency is $(2\pi)^{-1}$, as pointed out in Equation \eqref{eq: piston position},
the characteristic wavelength of the solution will be $L=2\pi$.
Different spatial refinement levels are considered, such
that $N$ cells are contained in a wavelength $L$,
for $N\in \{10,50,250\}$.
The amplitude of the piston's oscillations is set to $M=0.4$, which in 
the nonlinear case will lead
to a solution dominated by strong shocks, as studied by Inoue and Yano \cite{inoue_propagation_1993}.

Given a refinement level $N$ and a sponge layer length $\omega$,
we define
the error associated with an ABC
as
\begin{equation}
    \label{eq: Error ABC}
    E_{ABC}(N,\omega)=
    \max_{t\in [0,40\pi]}
    \frac{
    {
    \Vert
    u^{[N]}_{\text{ref}}(t)-u^{[N,\omega]}_{\text{ABC}}(t)
    \Vert
    }_{L^1(0,x_s)}
    }{\Vert u^{[N]}_{\text{ref}}(t)\Vert_{L^1(0,x_s)}},
\end{equation}
where we use a grid function norm over the cell-averaged velocity $u$,
and the reference solution is computed on a larger 
computational domain such that no wave interaction
takes place at the boundary.
Similarly, the error corresponding to the numerical discretization is approximated by
\begin{equation}
    \label{eq: Error Discretization}
    E_{num}(N)=
    \max_{t\in [0,40\pi]}
    \frac{
    {
    \Vert
    u^{[N]}_{\text{ref}}(t)-u^{[3000]}_{\text{ref}}(t)
    \Vert
    }_{L^1(0,x_s)}
    }{\Vert u^{[N]}_{\text{ref}}(t)\Vert_{L^1(0,x_s)}},
\end{equation}
where $u^{[3000]}_{\text{ref}}(\cdot)$ is undersampled by 
locally averaging the solution to match the grid of $u^{[N]}_{\text{ref}}(\cdot)$.

\subsection{Parameter configuration}
\label{sec: parameter config}

\begin{figure}[]
    \centering
    \begin{minipage}{0.48\textwidth}
        \centering
    \includegraphics[width=0.9\textwidth]{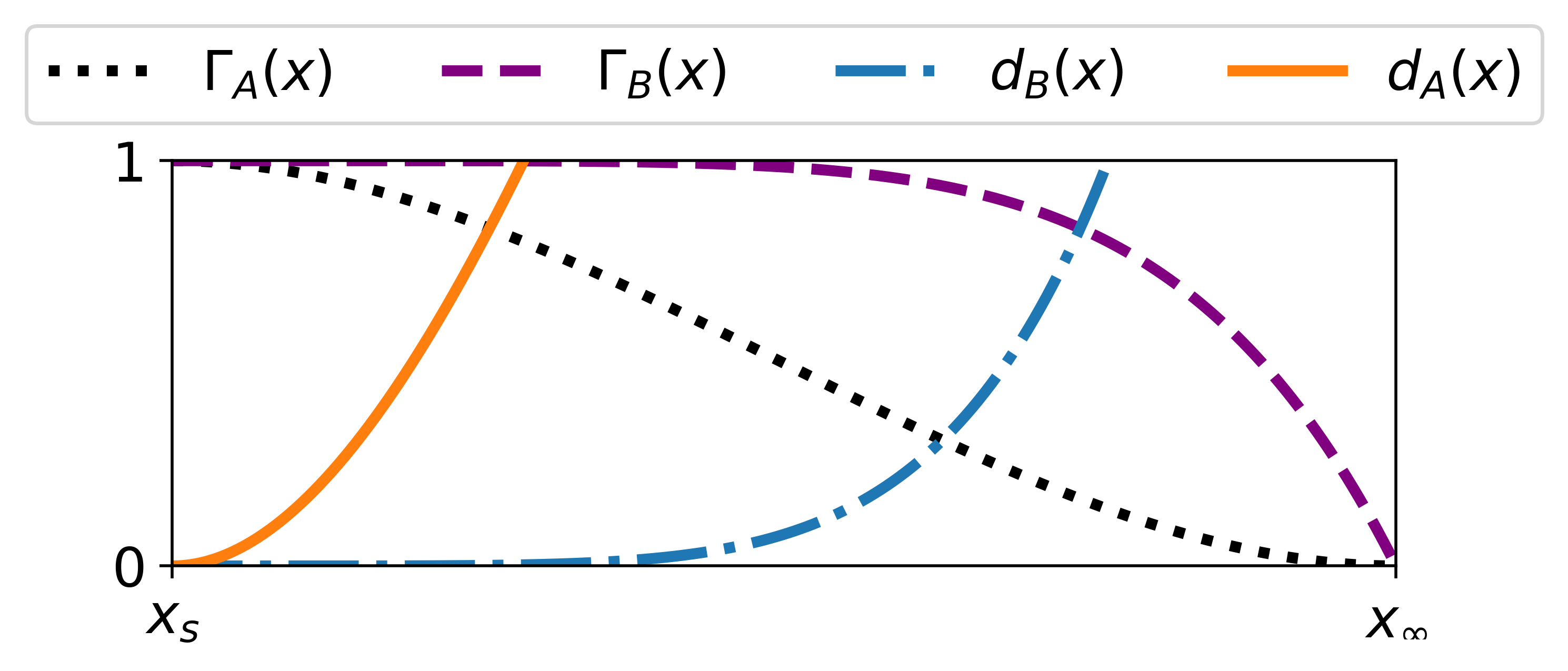}
    \caption{Weight functions $\Gamma_{(\cdot)}$ and damping functions $d_{(\cdot)}$, as defined in Equations \eqref{eq: Engsig-Karup weight function}, 
    \eqref{eq: Mayer weight function}, and \eqref{eq: definition d and s SDO} respectively}
    \label{fig: Shape slowing damping funcs}
    \end{minipage}%
    \hfill
    \begin{minipage}{0.5\textwidth}
        \centering
\includegraphics[width=\textwidth]{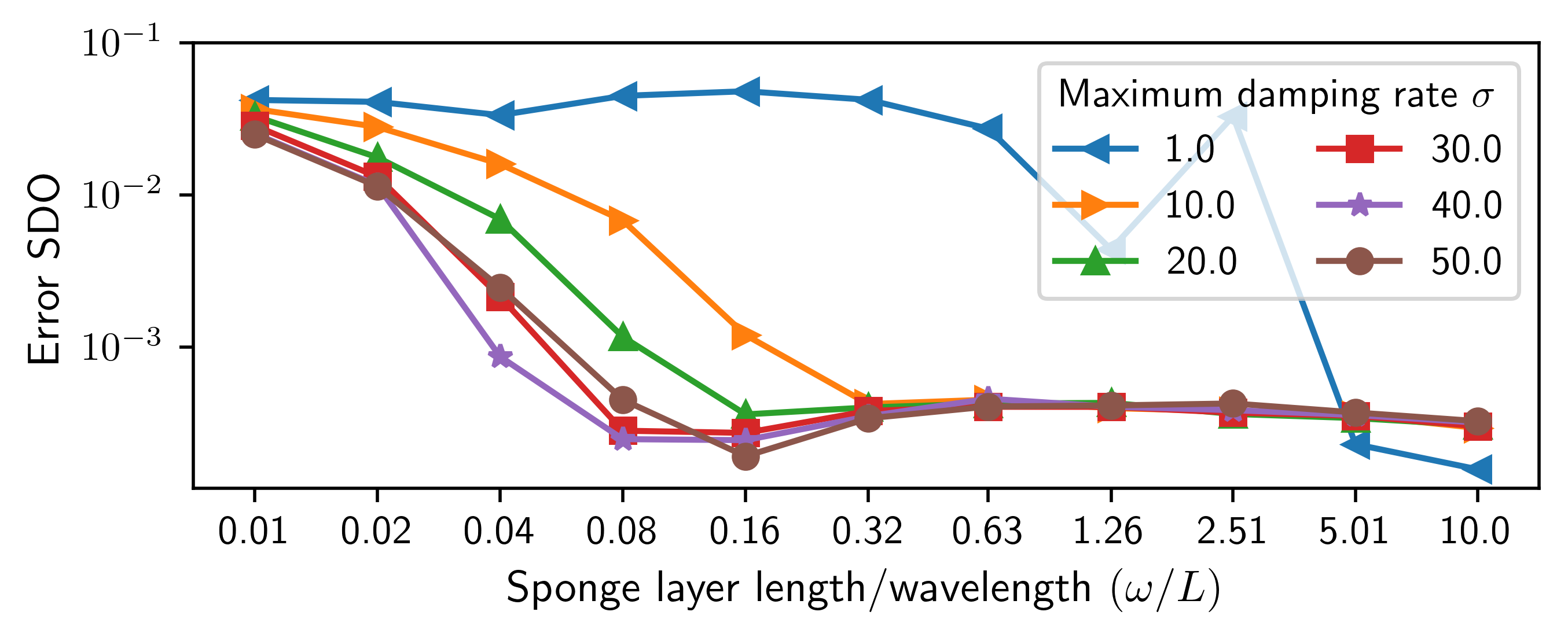}
        \caption{Relative error \eqref{eq: Error ABC} for the nonlinear equations \eqref{eq:conservative form lagrangian} with SDO ABC and different maximum damping rates $\sigma$ $(N=250)$}
        \label{fig:error different damping rates}
    \end{minipage}
\end{figure}

Before proceeding with a performance comparison between the
different ABCs, we tune the following parameters:
\begin{itemize}
    \item the slowing-down and damping functions $s$ and $d$, and the damping rate $\sigma$ \eqref{eq: definition d and s SDO}
     for the construction of the slowing-down and damping far-field operators 
    (SDO, S-SDO, NDO),
    \item a scalar weight function $\Gamma$ for the relaxation methods (RM, RM-M).
\end{itemize}
The functions $s$, $d$, and $\Gamma$ will be constructed using either $\Gamma_A$ \eqref{eq: Engsig-Karup weight function}
or  $\Gamma_B$ (with $b=1/2$) \eqref{eq: Mayer weight function}.
The slowing-down and damping functions constructed with $\Gamma_A$ ($\Gamma_B$) will be denoted by
 $s_A$ ($s_B$) and $d_A$ ($d_B$) respectively.
The shapes of these functions are illustrated in Figure \ref{fig: Shape slowing damping funcs}.

First, we evaluate different combinations of $s$ and $d$ for SDO with $\sigma=30$,
$N=250$, and consider different sponge layer lengths $\omega$.
The relative errors obtained with each combination are displayed in Table \ref{tab: Choosing slowing and damping SDO}.
For the linear equations \eqref{eq:linearized equations}, 
while the choice of $s$ and $d$ does not seem to be critical,
the best performance is attained with $s_A$ and $d_B$.
With this configuration, the slowing-down effect is greater than the damping at the early stage of the absorbing process.
Thus, the length of the outgoing waves is
reduced before they are damped and, consequently,
the drawbacks associated with large wavelengths, addressed in Section \ref{sec:Boundary treatment}, are alleviated.

\begin{table}[]
\centering
\begin{tabular}{ccccccccc}
\hline
\multirow{2}{*}{$\frac{\omega}{L}$} & \multicolumn{4}{c}{Linear equations \eqref{eq:linearized equations}} & \multicolumn{4}{c}{Nonlinear equations \eqref{eq:conservative form lagrangian}} \\ \cmidrule{2-9} 
                                     & $s_A$, $d_A$ & $s_A$, $d_B$ & $s_B$, $d_A$ & $s_B$, $d_B$ & $s_A$, $d_A$ & $s_A$, $d_B$ & $s_B$, $d_A$ & $s_B$, $d_B$ \\ \hline
$1/8$                                & 4.53e-04     & 4.46e-03     & 1.08e-03     & 2.74e-03     & 6.54e-04     & 5.63e-02     & 2.50e-04     & 2.66e-03     \\ \hline
$1/4$                                & 4.46e-04     & 5.29e-05     & 6.92e-04     & 3.62e-04     & 2.52e-01     & 6.87e-01     & 3.62e-04     & 4.53e-04     \\ \hline
$1/2$                                & 3.33e-04     & 5.11e-05     & 4.79e-04     & 1.66e-04     & 3.73e-01     & 3.86e-01     & 4.05e-04     & 3.96e-04     \\ \hline
$1$                                  & 2.03e-04     & 5.69e-05     & 2.44e-04     & 6.06e-05     & 4.10e-04     & 4.82e-01     & 4.13e-04     & 3.61e-04     \\ \hline
\end{tabular}
\caption{Error \eqref{eq: Error ABC} obtained with SDO and different
slowing-down and damping functions $(\sigma=30,\,N=250)$
\label{tab: Choosing slowing and damping SDO}}
\end{table}

\begin{figure}[htb]
    \centering
    \begin{minipage}{0.48\textwidth}
        \centering
\begin{tabular}{ccc}
\hline
\multirow{2}{*}{$\frac{\omega}{L}$} & \multicolumn{2}{c}{Nonlinear equations \eqref{eq:conservative form lagrangian}} \\ \cmidrule{2-3} 
                                     & $d_A$            & $d_B$            \\ \hline
$1/8$                                & 8.93e-03         & 5.74e-03         \\ \hline
$1/4$                                & 5.10e-03         & 3.48e-03         \\ \hline
$1/2$                                & 3.74e-03         & 2.80e-03         \\ \hline
$1$                                  & 3.21e-03         & 2.03e-03         \\ \hline
\end{tabular}
\captionof{table}{Relative error \eqref{eq: Error ABC} for the nonlinear equations \eqref{eq:conservative form lagrangian} with NDO ABC \ref{sec:NDO} and different damping functions $(\sigma=20,\,N=250)$ \label{tab: Different damping functions NDO}}
    \end{minipage}%
    \hfill
    \begin{minipage}{0.5\textwidth}
        \centering
\includegraphics[width=\textwidth]{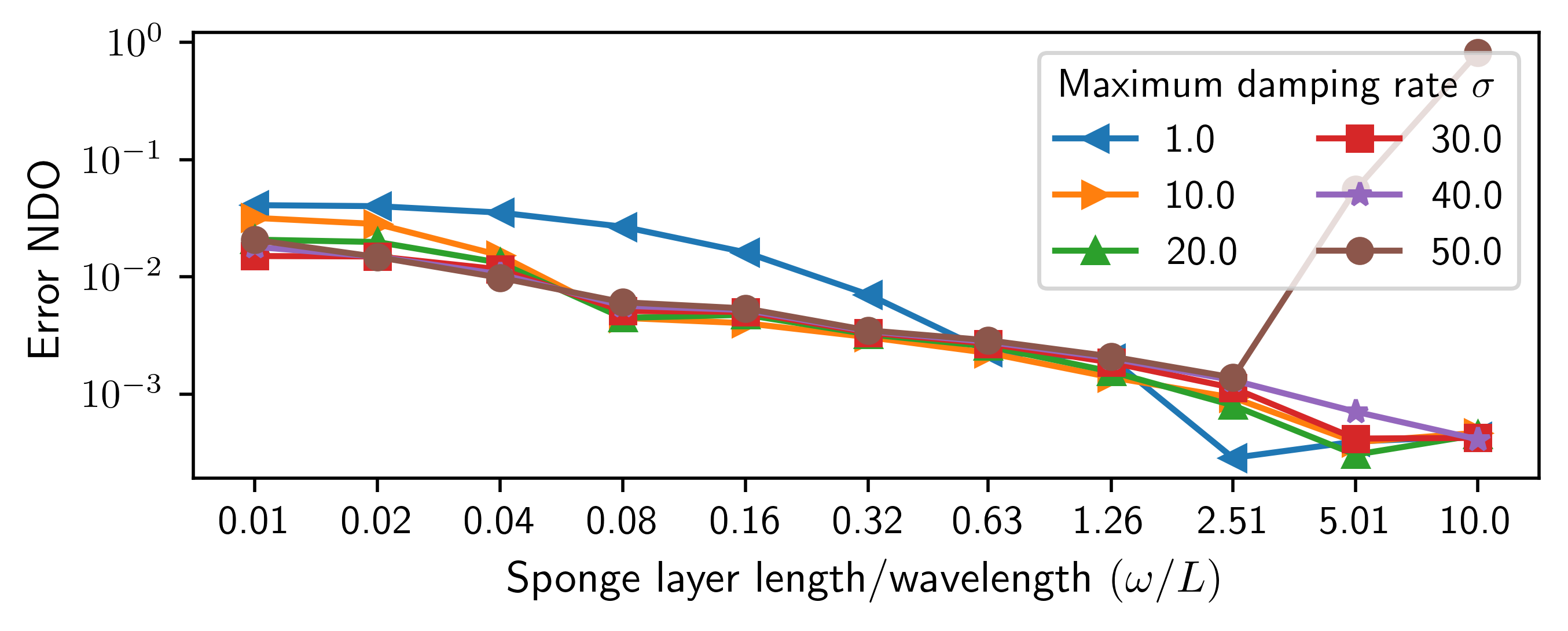}
        \caption{Relative error \eqref{eq: Error ABC} for the nonlinear equations \eqref{eq:conservative form lagrangian} with NDO ABC and different maximum damping rates $\sigma$ $(N=250)$}
        \label{fig:error different damping rates NDO}
    \end{minipage}
\end{figure}

In contrast with the linear scenario, the choice of slowing-down and damping functions 
drastically affects the behavior of the SDO ABC for the 
nonlinear equations \eqref{eq:conservative form lagrangian} 
(see Table \ref{tab: Choosing slowing and damping SDO}).
For this system, $s_B$ and $d_A$ yield the least error overall. 
This can be attributed to
the slowing-down operator's tendency to generate
entropy inside the sponge layer, illustrated in Appendix \ref{sec: Appendix stability slowing down},
which needs to be dissipated by a rapidly increasing damping rate.
The effect of using different maximum damping rates $\sigma$ with
SDO is depicted in Figure \ref{fig:error different damping rates}.
Since no significant improvement is obtained for larger values,
$\sigma=30$ is used for the SDO ABC in the following.
Proceeding analogously, we select $d_B$ and $\sigma=20$ for the NDO ABC, 
as suggested by the results displayed in Table \ref{tab: Different damping functions NDO} and Figure \ref{fig:error different damping rates NDO}.

Finally, the scalar weight function $\Gamma$ for the RM and RM-M techniques
is chosen similarly.
As can be seen in Table \ref{tab:RM different weight funcs},
the best performance is obtained with $\Gamma_B$.
Nevertheless, these techniques do not seem to be critically dependent on
the choice of $\Gamma$, even in the nonlinear setting.

Although some of these parameters could be optimized through reflection-coefficient-based strategies \cite{peric2018analytical},
in this work we tune them based on the error \eqref{eq: Error ABC} only to ensure a fair and consistent performance comparison
among the different ABCs.

\begin{table}[]
\centering
\begin{tabular}{ccccccccc}
\hline
\multirow{3}{*}{$\frac{\omega}{L}$} & \multicolumn{4}{c}{Linear equations \eqref{eq:linearized equations}} & \multicolumn{4}{c}{Nonlinear equations \eqref{eq:conservative form lagrangian}} \\ \cmidrule{2-9} 
                                     & \multicolumn{2}{c}{RM}               & \multicolumn{2}{c}{RM-M}            & \multicolumn{2}{c}{RM}                  & \multicolumn{2}{c}{RM-M}             \\ \cmidrule{2-9} 
                                     & $\Gamma_A$       & $\Gamma_B$ & $\Gamma_A$ & $\Gamma_B$ & $\Gamma_A$          & $\Gamma_B$ & $\Gamma_A$ & $\Gamma_B$ \\ \hline
$1/8$                                & 3.30e-03         & 2.16e-03   & 8.64e-04   & 1.36e-03   & 8.06e-04            & 9.69e-04   & 1.24e-03   & 5.18e-04   \\ \hline
$1/4$                                & 1.94e-03         & 7.75e-04   & 6.54e-04   & 5.04e-04   & 1.06e-03            & 1.09e-03   & 3.78e-04   & 2.42e-04   \\ \hline
$1/2$                                & 8.78e-04         & 2.31e-04   & 4.40e-04   & 2.23e-04   & 1.10e-03            & 2.13e-03   & 2.02e-03   & 4.61e-04   \\ \hline
$1$                                  & 3.05e-04         & 6.76e-05   & 3.05e-04   & 8.71e-05   & 1.09e-03            & 9.81e-04   & 4.37e-04   & 4.13e-04   \\ \hline
\end{tabular}
\caption{Error \eqref{eq: Error ABC} obtained with the scalar \ref{sec:RM} and matrix-valued \ref{sec: relaxation matrix valued} RMs and different weight functions $(N=250)$\label{tab:RM different weight funcs}}
\end{table}

\subsection{Performance comparison}
\label{sec: performance comparison}

\begin{figure}[htb]
    \centering
    \begin{minipage}{0.49\textwidth}
        \centering
    \includegraphics[width=.9\textwidth]{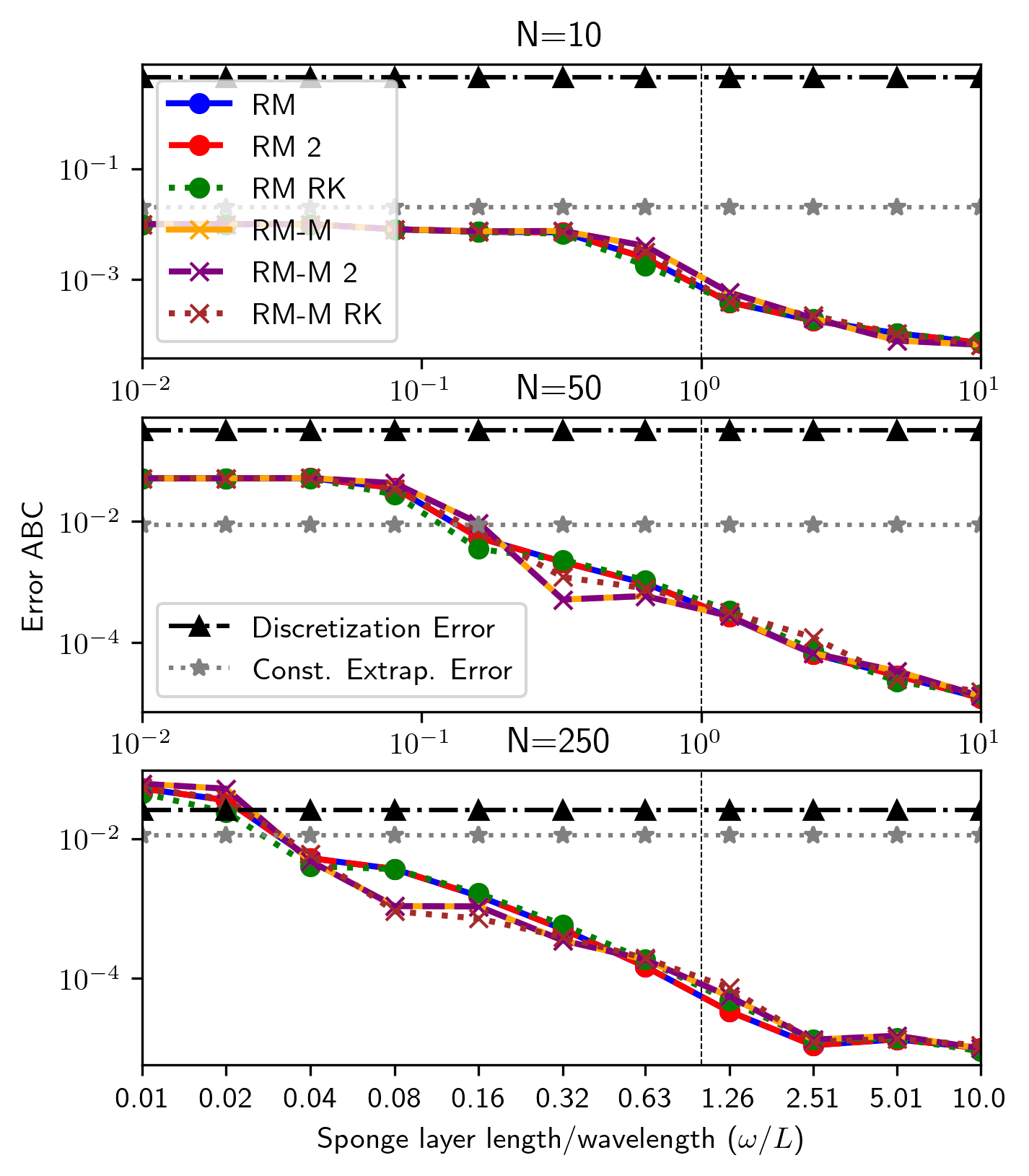}
    \caption[width=.9\textwidth]{Comparison of relaxation methods for the linear equations \eqref{eq:linearized equations}}
    \label{fig: comparison RM linear}
    \end{minipage}%
    \hfill
    \begin{minipage}{0.49\textwidth}
        \centering
\includegraphics[width=.9\textwidth]{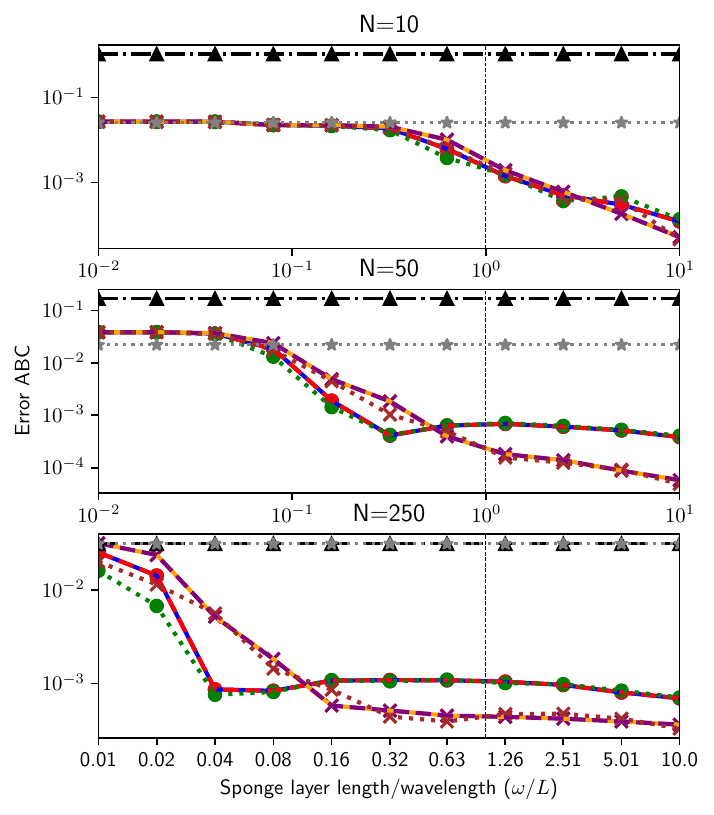}
        \caption[width=.9\textwidth]{Comparison of relaxation methods for the nonlinear equations \eqref{eq:conservative form lagrangian}}
        \label{fig: comparison RM nonlinear}
    \end{minipage}
\end{figure}

Let us compare the behavior of the different 
relaxation methods discussed in  Sections \ref{sec:RM}, \ref{sec: relaxation matrix valued}, and \ref{sec: Strang}.
We refer to Algorithms \ref{alg: relaxation scalar} and \ref{alg: relaxation matrix}
as RM and RM-M respectively, and to their use described
in Algorithm \ref{alg: Strang splitting} as RM 2 and RM-M 2.
For the latter, the weight function $\Gamma_B^2$ is used,
according to the observations from Section \ref{sec: Strang}.
Moreover, we consider the application of these techniques after each stage of the
Runge-Kutta time integration, as proposed by Wasistho et al. \cite{wasistho_simulation_1997},
and refer to them as RM RK and RM-M RK.
The results are displayed in Figures \ref{fig: comparison RM linear} and
 \ref{fig: comparison RM nonlinear} for the linear \eqref{eq:linearized equations} 
and nonlinear \eqref{eq:conservative form lagrangian}
equations respectively.

To isolate the effect of the ABCs, the error associated with the numerical discretization
\eqref{eq: Error Discretization} is displayed in black dot-dashed lines with triangle markers
in all the figures of this section.
Furthermore, for the sake of comparison, the error associated with the use of constant extrapolation at the last 
ghost cell without a sponge layer (a simple and common boundary treatment)
is shown in grey dotted lines with star markers.

For the linear equations, there is no perceivable difference between RM and RM-M, and their variants.
This can be attributed to the fact that due to the piston motion, all the waves that enter 
the sponge layer are right-going, and no interaction between them takes place in the linear setting.
Therefore, there is no difference between absorbing waves selectively or not.
Conversely, for the nonlinear equations, the only noticeable difference is between the scalar (RM)
and directional variants (RM-M) of the relaxation method when larger portions of 
a characteristic wavelength $L$ are contained
in the sponge layer.

\begin{figure}[htb]
    \centering
    \begin{minipage}{0.49\textwidth}
       \centering
        \includegraphics[width=.9\textwidth]{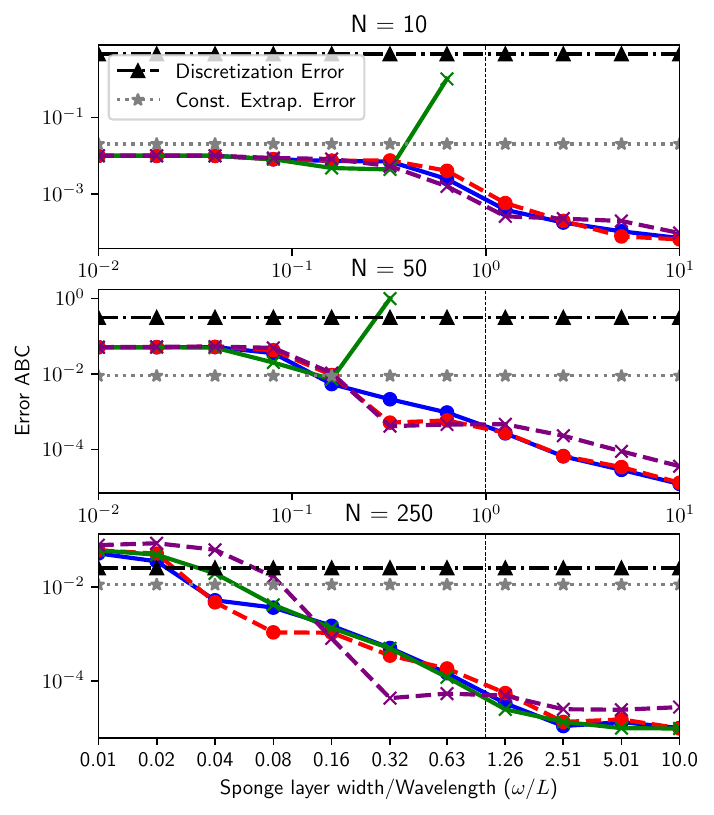}
        \caption[width=.9\textwidth]{Comparison of different ABCs for the linear equations \eqref{eq:linearized equations}}
        \label{fig:error comparison all methods linear}
    \end{minipage}%
    \hfill
    \begin{minipage}{0.49\textwidth}
        \centering
        \includegraphics[width=.9\textwidth]{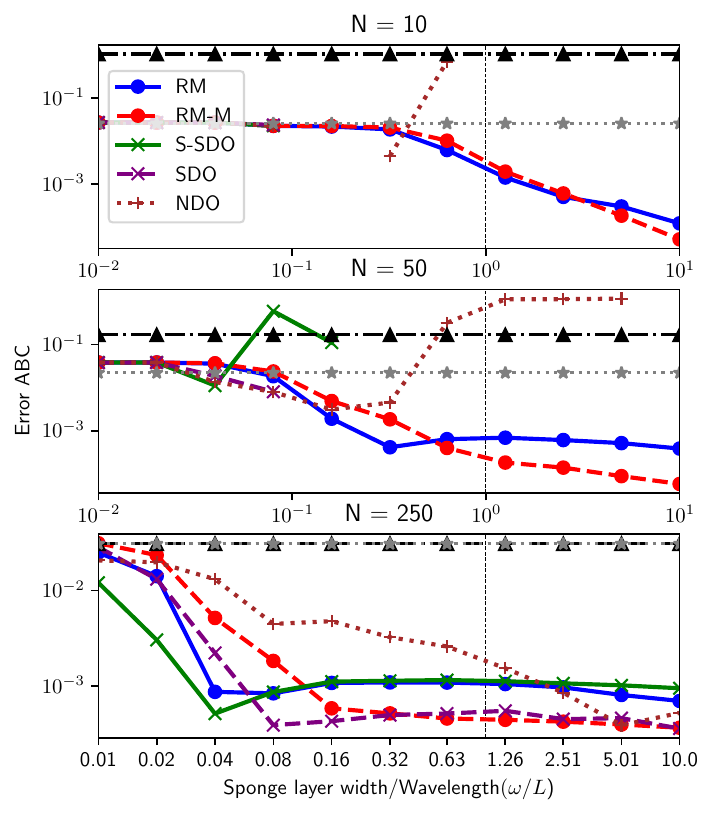}
        \caption[width=.9\textwidth]{Comparison of different ABCs for the nonlinear equations
        \eqref{eq:conservative form lagrangian}
        \label{fig:error comparison all methods nonlinear}}
    \end{minipage}
\end{figure}

Now we assess the perfomance of all the ABCs addressed in Sections \ref{sec:Boundary treatment}
and \ref{sec: Extensions}; the results 
are displayed in Figures \ref{fig:error comparison all methods linear} and
\ref{fig:error comparison all methods nonlinear}. 
The methods based on far-field operators
(SDO, S-SDO, NDO)
fail to provide a satisfactory error and the numerical solution diverges
for the coarsest grids.
This suggests that the discretization is not stable with the time step size prescribed by the CFL condition
with the chosen Courant number.
Indeed, we can see that the $L^2$ norm of $\bfD(x)$ in Equation \eqref{eq: linear SDO} 
is $\displaystyle \max_i|d_i s_i \lambda_i|$, 
which is larger than time step size restriction
for a coarse grid when the maximum damping parameter $\sigma$ is large.
As seen in Figure \ref{fig:error different damping rates},
a large maximum damping rate might be necessary to attain a prescribed error
for a fixed sponge layer length.
Since it is convenient to use an explicit scheme for the hyperbolic part of the 
problem, a way to deal with this artificial stiffness, 
not explored in this work, is to use
an ImEx scheme \cite{pareschi_implicit-explicit_2005} 
to evolve the semidiscretization.
On the other hand, the RM and its extension RM-M, introduced herein, 
show to be the most robust ABCs for both the linear and nonlinear equations.

Once the grid (and thus the time step size) is refined, and the stability issues
are overcome, we can see that the SDO ABC attains its minimum error for a sponge
layer length of around $0.08L$.
It is followed by RM-M and NDO which attain their minimum error 
for sponge layer lengths of $0.6L$ and $5L$ respectively.
While the scalar methods S-SDO and RM attain a low error for a sponge layer length of $0.04L$,
they still do not reach the absorption of the other methods, even
for a sponge layer length of $10L$.
Indeed, despite the fact that all the ABCs eventually
attain an error that is approximately $1.5$ orders of magnitude smaller than
the numerical discretization error,
the directional methods (SDO, RM-M) further improve upon
the scalar ones (S-SDO, RM) by about half an order of magnitude for the nonlinear equations.
This can be explained by the fact that, in the nonlinear setting,
 the first two shocks of the wave train
are periodically merged in the far-field \cite{inoue_propagation_1993},
and such interactions might cause information to propagate back into the
computational domain \cite{hedstrom_nonreflecting_1979}.
Consequently, as the sponge layer size increases,
methods that selectively damp a characteristic field, like SDO, NDO, and RM-M, allow this
physical information to recirculate, while scalar techniques, like RM and S-SDO,
fail to capture this behavior.
Even though a sponge layer length much larger than the characteristic wavelength $L$ is not
generally used in practice, it is useful in Figures \ref{fig:error comparison all methods linear}
and \ref{fig:error comparison all methods nonlinear} to illustrate that a larger sponge layer
does not compensate for the lack of selectivity of scalar methods.

Although unphysical values (e.g. negative pressure) are not observed with
any method, it is worth mentioning that,
being a convex combination of admissible states, the RM is the only method that is a priori guaranteed to
preserve physical admissibility of the solution (as long as the underlying numerical scheme does so). 
For close-to-vacuum states, further care might be required for the other methods to maintain
physical values.

Finally an execution time comparison between the different ABCs
for the nonlinear equations \eqref{eq:conservative form lagrangian}
is displayed in Table \ref{tab: execution time}.
We consider the best wall clock time obtained for each method after 7 runs 
on a
MacBook Pro 2023 with a 10-Core CPU Apple M2 Pro chip.
As expected, the scalar methods, RM and S-SDO, are the fastest
ones, closely followed by the directional methods
RM-M and SDO.
The NDO technique is, by a large margin, the most
computationally expensive due to the treatment of its nonlocal integral term
(see Appendix \ref{apen: Discretization NDO}).

\begin{table}
    \centering
    \begin{tabular}{cccccc}
    \hline
        N & RM & RM-M & S-SDO & SDO & NDO \\
    \hline
        10 & 0.12 & 0.14 & 0.12 & 0.14 & 0.15 \\
    \hline
        50 & 0.56 & 0.79 & 0.57 & 0.66 & 1.78 \\
    \hline
        250 & 6.91 & 11.72 & 8.12 & 11.40 & 40.96 \\
    \hline
    \end{tabular}
    \caption{Execution time (s) comparison of ABCs with different refinement levels $N$ for the nonlinear equations
         \label{tab: execution time} with
        a sponge layer of length $\omega=10L$}
    \label{tab:execution_time}
\end{table}

\section{Conclusion}
Methods from two families of absorbing boundary conditions were
studied, and their performance was tested
for the piston problem of gas dynamics.
The relaxation method with a scalar-valued weight function is the least expensive and simplest to implement.
Karni's slowing-down and damping operators are able to selectively absorb waves from a 
specific characteristic field, allowing physical waves to circulate back into the computational
domain.
However, this technique requires the addition of an artificial source term that might
become stiff depending on the prescribed damping rate, leading to additional time step size
restrictions to maintain stability with an explicit time integrator.
Among the extensions for both techniques presented in this work, 
the relaxation method with a matrix-valued weight function is an attractive alternative,
as it preserves the remarkable robustness of the relaxation method (even in the nonlinear setting),
while selectively damping outgoing waves.
We expect this property to be more advantageous in a case where more interactions take place 
inside the sponge layer, for instance, when the sponge layer is used as a wavemaker 
(i.e. nonconstant far-field state), which is a common application for the relaxation method.
Furthermore, sophisticated scalar weight functions (e.g. that depend on the conserved quantities),
that have already been developed for the
relaxation method, can be used with this extension in a straightforward manner.
Due to its higher computational cost and lack of robustness for coarse grids, the newly proposed nonlinear damping operator
is not competitive with existing ABCs considered in this work.
However, this technique could warrant further investigation
in scenarios where a higher-order representation of the solution is available inside each cell,
since the integral damping term could leverage this additional information.
While this work focused on one-dimensional problems, the extension of the matrix-valued
relaxation method to multiple spatial dimensions could be achieved for
rotationally invariant systems following the approach in \cite{karni_far-field_1996};
this is the subject of future work.

\section*{Acknowledgements}

I would like to express my deepest gratitude 
to Prof. David Ketcheson  and Prof. Giovanni Russo 
for their invaluable advice and 
support throughout the development of this work.
This work has been funded by King Abdullah University of Science
and Technology.

\begin{appendices}
\section{Reflection coefficients}
\label{sec: reflection analysis}
Based on the relations among the different ABCs presented in Section \ref{sec: Extensions},
it is possible to provide analytical estimates of the absorbing properties of 
the SDO and RM techniques in a unified framework.
For this, we restrict ourselves to linear systems of the form \eqref{eq:conservation law linear}
and begin by considering the SDO \ref{sec: SDO} with no
slowing down, i.e., $s_i(x)\equiv 1$ for all $i$.
From \eqref{eq: explicit solution  SDO characteristic no slowing}
we can see that, if  $\int_{x_s}^{x}d_i(z) \, dz \to +\infty$ as $x \to x_\infty$,
the  right-going information will
be completely absorbed inside the sponge layer. 
In practice, however, $\max(d_i)=\sigma_i$ is chosen to be large but finite to avoid numerical stiffness.
Thus, reflection of waves will occur at $x_\infty<+\infty$.

In a similar manner as done in \cite{peric2018analytical} by Peri\'c and
Abdel-Maksoud for scalar-type damping operators, 
we 
assume that the boundary at $x_\infty$ is perfectly reflecting
to estimate the reflection coefficient of the sponge layer associated with SDO.
We consider an initial condition that is compaclty supported in the computational
domain $[0,x_s]$, which has a characteristic wavelength $L$,
and is centered at $\mathring{x}$ ($\mathring{x}+L/2<x_s$).
Moreover, the incoming wave is is purely right-going, i.e.,
\begin{equation}
    {\bfq}(0,x) =  \bfq^{\text{inc}}(x)= f(x)\chi_{[\mathring{x}-L/2,\mathring{x}+L/2]}(x) \sum_{\lambda_i>0} a_i \bfr_i,
\end{equation}
for some coefficients $a_i$.
In characteristic variables, we allow the solution to evolve up to a point
where the information in each characteristic field has left
the sponge layer, i.e., 
$$t_\infty = \min_{t>0}\max_{\lambda_i>0}(\mathring{x}-L/2+\lambda_i t)>x_\infty.$$
Then, for each $x>x_\infty$, and for each characteristic field $i$,
we have
\begin{align}
    w_i(t_\infty,x) &= a_i f(x-x_i)\chi_{[x_i-L/2,x_i+L/2]}(x)
    \exp\left(-\int_{x-\lambda_i t_\infty}^{x} d_i(z) \, dz\right)\\
    & = a_i f(x-x_i)\chi_{[x_i-L/2,x_i+L/2]}(x)
    \exp\left(-\int_{x_s}^{x_\infty} d_i(z) \, dz\right),
\end{align}
where $x_i=\mathring{x}+\lambda_i t_\infty$.

Let us assume that the eigenstructure of the system is symmetric, like
in the case of the linearized Euler equations in Lagrangian coordinates,
presented in Section \ref{sec:Euler equations} below.
Consequently, the information from each characteristic field $i$ is reflected back
through its reciprocal characteristic field $j(i)$ such that $\lambda_j=-\lambda_i$,
after interacting with the perfectly reflecting boundary at $x_\infty$.
Thus, the total mass (of the $k$-th conserved variable) of the reflected wave
is given by
\begin{align}
    \label{eq: reflected mass}
    M_{\text{ref}} &=  || \bfq^{\text{ref}}_k(t_\infty,\cdot)||_{L^1}\\
    &= \sum_{\lambda_j<0}||f(x-x_i)\chi_{[x_i-L/2,x_i+L/2]}(x)||_{L^1} |a_j| |\bfr_{jk}|  \exp\left(-\int_{x_s}^{x_\infty} d_i(z) \, dz\right)\\
    &= ||f(\cdot)\chi_{[\mathring{x}-L/2,\mathring{x}+L/2]}(\cdot)||_{L^1}\sum_{\lambda_j<0} |a_j| |\bfr_{jk}|  \exp\left(-\int_{x_s}^{x_\infty} d_i(z) \, dz\right),
\end{align}
while the total mass of the incoming wave is 
\begin{align}
    \label{eq: incoming mass}
    M_{\text{inc}} =||f(\cdot)\chi_{[\mathring{x}-L/2,\mathring{x}+L/2]}(\cdot)||_{L^1} \sum_{\lambda_i>0}|a_i| |\bfr_{ik}|.
\end{align}
Then, we define the reflection coefficient associated with the damping operator as
\begin{equation}
    C_R^{\text{SDO}} = \frac{M_{\text{ref}}}{M_{\text{inc}}} = \frac{\sum_{\lambda_j<0}|a_j| |\bfr_{jk}| \exp\left(-\int_{x_s}^{x_\infty} d_i(z) \, dz\right)}{\sum_{\lambda_i>0}|a_i| |\bfr_{ik}|},
\end{equation}
which, due to the symmetry of the eigenstructure, can be rewritten as
\begin{equation}
    C_R^{\text{SDO}} = \frac{\sum_{\lambda_i>0}|a_i| |\bfr_{ik}| \exp\left(-\int_{x_s}^{x_\infty} d_i(z) \, dz\right) }{\sum_{\lambda_i>0}|a_i|  |\bfr_{ik}|}.
\end{equation}

Following a similar reasoning, we can compute the reflection coefficient of the
scalar damping operator \eqref{eq: scalar SDO}.
In this case, 
after the wave is reflected at $x_\infty$, the  left-going information
is not only propagated back by left-going characteristics $\lambda_j<0$, but
also damped at the rate 
$$
\exp \left(-\mathrm{sgn}(\lambda_j)|\lambda_j^{-1}|\int_{x-\lambda_j t}^{x} d(z) \, dz\right)= 
\exp \left(\lambda_j^{-1}\int_{x-\lambda_j t}^{x} d(z) \, dz\right).
$$
Hence, the total mass of the reflected wave is
\begin{align}
    \label{eq: reflected mass sSDO}
    M_{\text{ref}} 
    = ||f(\cdot)\chi_{[\mathring{x}-L/2,\mathring{x}+L/2]}(\cdot)||_{L^1}\sum_{\lambda_j<0} |a_j| |\bfr_{jk}| 
    \exp\left((\lambda_j^{-1}-\lambda_i^{-1})\int_{x_s}^{x_\infty} d(z)\, dz\right),
\end{align}
which, due to the symmetry of the eigenstructure, gives us the reflection coefficient
\begin{equation}
    C_R^{\text{S-SDO}} =  \frac{
    \sum_{\lambda_i>0}|a_i| |\bfr_{ik}| \exp\left(-2\lambda_i^{-1}\int_{x_s}^{x_\infty} d(z) \, dz\right)
     }{\sum_{\lambda_i>0}|a_i|  |\bfr_{ik}|}.
\end{equation}
We can see that the reflection coefficient of S-SDO will be smaller
or larger than that of SDO depending on the magnitude of the eigenvalues $\lambda_i$.

Using this setup, we can also compute the reflection coefficient for the matrix-valued
relaxation method RM-M.
We assume that left-going waves are unperturbed ($\Gamma_j\equiv 1$ for $j<0$),
and consider the damping function \eqref{eq: log damping rate},
which necessarily depends on the discretization through $\Delta t$.
Then, proceeding as above, the reflection coefficient of this method
is given by 
\begin{align}
    C_R^{\text{RM-M}} 
    & = \frac{\sum_{\lambda_i>0}|a_i||\bfr_{ik}|
    \exp\left(\frac{1}{\Delta t \lambda_i}\int_{x_s}^{x_\infty} \ln\circ \Gamma_i(z)\, dz\right) }{\sum_{\lambda_i>0}|a_i|  |\bfr_{ik}|}\\
    & \approx \frac{\sum_{\lambda_i>0}|a_i||\bfr_{ik}|
    \exp\left(\frac{x_\infty-x_s}{C\Delta x}\int_{0}^{1} \ln\circ \Gamma_i\circ\phi^{-1}(z)\, dz\right) }{\sum_{\lambda_i>0}|a_i|  |\bfr_{ik}|}.
\end{align}
Here, we have used the fact that $\phi$ is invertible
in $[x_s,x_\infty]$ \eqref{eq:scaling weight functions}, and assumed the CFL condition
$C=\Delta t \max_{i,p} |\lambda_p|/\Delta x$.
The factor $(x_\infty-x_s)/(C\Delta x)$ indicates that the effectiveness of the RM-M method
depends not only on the choice of $\Gamma_i$, but also on the number of cells in the sponge layer. 
Similarly as before, we can compute the reflection coefficient for the scalar relaxation method RM,
which is given by
\begin{align}
    C_R^{\text{RM}} \approx \frac{\sum_{\lambda_i>0}|a_i||\bfr_{ik}|
    \exp\left(2\frac{x_\infty-x_s}{C\Delta x}\int_{0}^{1} \ln\circ \Gamma_i\circ\phi^{-1}(z)\, dz\right) }{\sum_{\lambda_i>0}|a_i|  |\bfr_{ik}|}.
\end{align}
As can be expected from our assumptions, if information from all the characteristic fields
is damped with the scalar method RM, the reflection coefficients will be smaller
than those of the directional methods RM-M.

\subsection{Comparsion with numerical results}
To validate the theoretical estimates of the reflection coefficients
derived in the previous section,
we consider the linearized equations \eqref{eq:linearized equations}.
For this, we take $x_s=20\pi$, different values of $x_\infty$.
The initial condition is given by the purely right-going pulse
\begin{align}
    \bfq(x,0) = \sin(x-10\pi)\chi_{[9\pi,11\pi]}(x)\bfr_3,
\end{align}
and the solution is evolved up to $t_F=50s$.
Reflecting boundary conditions are set at $x_\infty$ by changing
the sign of the velocity component of the solution in the rightmost ghost cell.
The numerical reflection coefficient
is estimated as the ratio
\begin{align}
    C^{\text{num}}_R = \frac{\max_{x\in[0,x_s]} u(t_F,x)-\min_{x\in[0,x_s]} u(t_F,x)}{\max_{x\in[0,x_s]} u(0,x)-\min_{x\in[0,x_s]} u(0,x)}.
\end{align}

\begin{figure}[]
    \centering
    \begin{minipage}{0.48\textwidth}
        \centering
    \includegraphics[width=\textwidth]{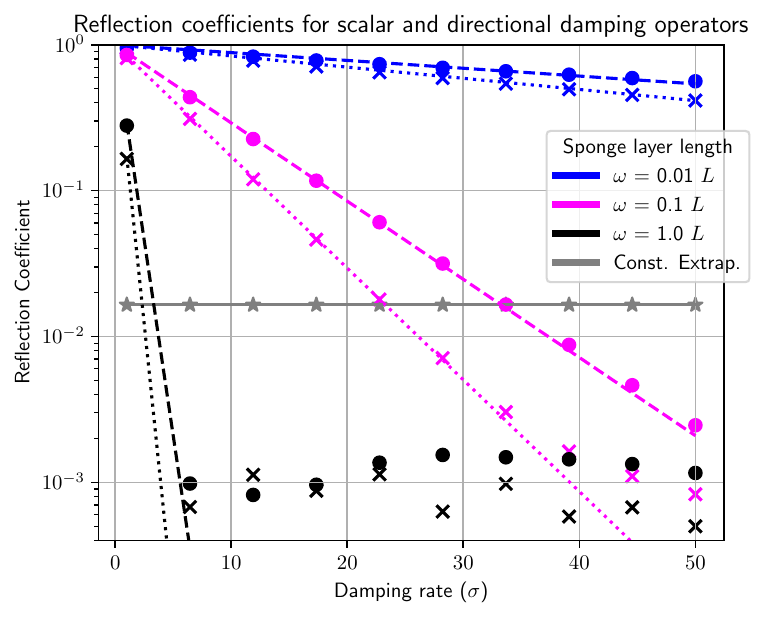}
      \caption{Comparsion of numerical (scattered points) and theoretical (lines) reflection coefficients 
    for the linear equations \eqref{eq:linearized equations} with 
    SDO (circle markers and dashed lines) and S-SDO ABCs (cross markers and dotted lines) respectively.}
    \label{fig: reflection SDO}
    \end{minipage}%
    \hfill
    \begin{minipage}{0.48\textwidth}
        \centering
\includegraphics[width=\textwidth]{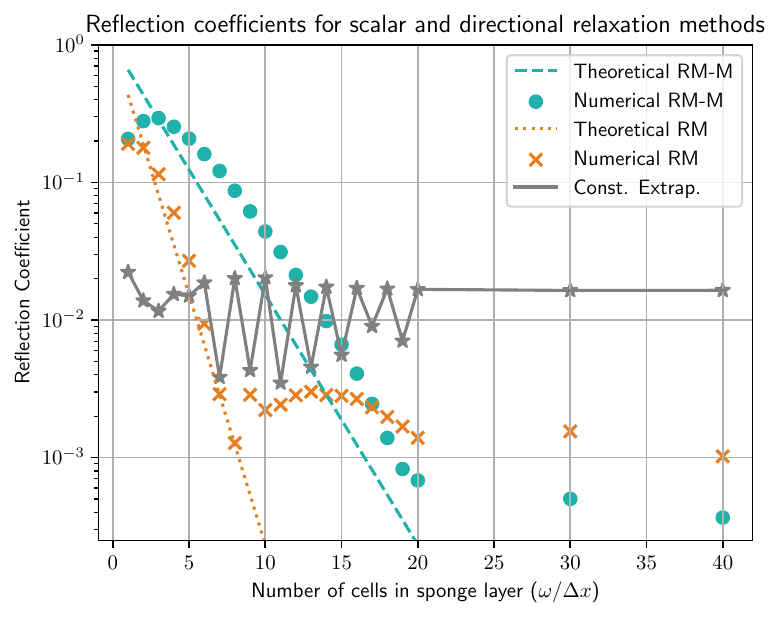}
        \caption{Comparsion of numerical (scattered points) and theoretical (lines) reflection coefficients 
    for the linear equations \eqref{eq:linearized equations} with 
    RM-M (circle markers and dashed lines) and RM ABCs (cross markers and dotted lines) respectively.}
        \label{fig: reflection RMM}
    \end{minipage}
\end{figure}

For the far-field operators SDO and S-SDO, we consider sponge layers of lengths $\omega \in \{0.01L,0.1L, L\}$ for
$L=2\pi$, and  $s_i=1$ for all $i\in\{1,2,3\}$.
For S-SDO we take $d_i$  as in \eqref{eq: definition d and s SDO}, with 
$\Gamma_d = \Gamma_B$ \eqref{eq: Mayer weight function}.
We proceed analoglously for SDO, but with $d_i\equiv 0$ for $i=1,2$.
Different damping rates $\sigma_i$ are considered to evaluate their influence on the reflection coefficient.
Similarly, for the relaxation methods RM-M and RM, we take the weight function $\Gamma=\Gamma_B$ as in \eqref{eq: Mayer weight function}.
Since for these methods the reflection coefficients depend only on the number of cells in the sponge layer, 
we consider a fixed sponge layer length $\omega=L$ and different refinement levels $\Delta x$.
For all the tests, good agreement between the numerical and theoretical reflection coefficients is observed
until saturation is attained,
as shown in Figures \ref{fig: reflection SDO} and \ref{fig: reflection RMM}
and discussed below.

While the theoretical estimates could inform the choice of parameters for the ABCs,
it is important to note that the conditions under which they were derived
are more often than not violated in practice.
For instance, in a realistic simulation, additional care is taken to absorb
waves at the end of the sponge layer $x_\infty$, by imposing 
the far-field state $\overline{\bfq}$ in the ghost cells or using constant extrapolation
instead of perfectly reflecting boundary conditions.
A numerical estimate of the reflection coefficient associated with the use of constant extrapolation
(with no sponge layer) is displayed in Figures \ref{fig: reflection SDO} and \ref{fig: reflection RMM}
(note that different resolutions for the numerical scheme are used in Figure \ref{fig: reflection RMM}, thus
constant extrapolation has a slightly different effect for each $\Delta x$).
Moreover, waves are not always purely outgoing when they reach the sponge layer,
which can be even more impactful in the nonlinear case, where characteristic fields
are continuously interacting.

Finally, the theoretical estimates presented here for the SDO and RM techniques are
insensitive to the numerical discretization.
Most importantly, we derive these estimates under the severe assumption
that reflections only occur at the end of the sponge layer $x_\infty$,
while in practice, reflections continuously take place inside the sponge layer.
These discrepancies become evident for the SDO ABC in Figure \ref{fig: reflection SDO}, as the 
numerical reflection coefficients suddenly saturate as the damping rates $\sigma$ increase,
and analogously for the relaxation methods in Figure \ref{fig: reflection RMM} as the grid is refined.
These effects must be taken into account when using the theoretical estimates
to guide the choice and optimization of parameters for the ABCs, see e.g. \cite{peric2018analytical,peric2020reducing}.

\label{sec: Numerical reflection coefficients}

\section{Construction of damping operators}
\label{apen: Construction of operators}
The construction of the damping operators and the matrix-valued weight function of
Sections \ref{sec: SDO}, \ref{sec:NDO}, and \ref{sec: relaxation matrix valued}, relies
on the computation of a matrix $\bfM=\bfR \boldsymbol{\Sigma}\bfR^{-1}$, where
$\boldsymbol{\Sigma}=\text{diag}(\sigma_i)$, and $\bfR$ contains
the right eigenvectors of the flux Jacobian matrix as columns.
For the linear equations \eqref{eq:linearized equations}
we have
\begin{equation}
 \bfM=
 \left( \begin{array}{ccc}
\frac{1}{\gamma}\left(\frac{p_1+\sigma_3}{2}+(\gamma-1)\sigma_2\right) & \frac{1}{2\gamma}(\sigma_1-\sigma_3) & \frac{\gamma-1}{\gamma^2}\left(\frac{-\sigma_1-\sigma_3}{2}+\sigma_2\right) \\
\frac{1}{2}(\sigma_1-\sigma_3) & \frac{1}{2}(\sigma_1+\sigma_3) & \frac{\gamma-1}{2\gamma}(\sigma_3-\sigma_1) \\
-\frac{1}{2}(\sigma_1+\sigma_3)+\sigma_2 & \frac{1}{2}(\sigma_3-\sigma_1) & \frac{1}{\gamma}\left(\frac{\gamma-1}{2}(\sigma_1+\sigma_3)+\sigma_2\right)
\end{array} \right),
\end{equation}
while for the nonlinear equations \eqref{eq:conservative form lagrangian}
we get
\begin{align}
\bfM=
\begin{pmatrix}
\frac{1}{2\gamma}(2\gamma \sigma_2+C_1) & \frac{1}{2p\gamma}(\sqrt{\gamma p V} C_2+u C_3 C_1) & \frac{1}{2p\gamma}(-C_1 C_3) \\
\frac{C_2 p}{C_4} & \frac{\sqrt{\gamma} \sqrt{p V}(\sigma_1+\sigma_3)+u C_3 C_2}{C_4} & \frac{-C_2 C_3}{C_4} \\
\frac{-pV C_1+C_2\sqrt{p V}\sqrt{\gamma}u}{2V\gamma} & \frac{\sqrt{p V}uC_1+\sqrt{\gamma}(C_3u^2-pV)C_2}{C_4\sqrt{\gamma}} & \frac{-\sqrt{p V}C_1-\sqrt{\gamma}C_3C_2u+\gamma \sqrt{p V}(\sigma_1+\sigma_3)}{C_4\sqrt{\gamma}}
\end{pmatrix},
\end{align}
where $C_1 = \sigma_1-2\sigma_2+\sigma_3$, $C_2 = \sigma_1-\sigma_3$,
$C_3 = \gamma-1$, and $C_4 = \sqrt{\gamma p V}$.

\section{Discretization of integral source term with nonconservative product}
\label{apen: Discretization NDO}
 Let us consider the PDE
 \begin{align}
    \label{eq: nonconservative integral}
     \partial_t \bfq(t,x) + \partial_x\bff(\bfq(t,x) ) = \int_x^{x_\infty} \bfM(z,\bfq) \partial_z \bfq(t,z)\,dz,
 \end{align}
 under the Ansatz that the solution is piecewise affine
 inside the cells $\mathcal{C}_i$ and not necessarily
 continuous at the interfaces, e.g.,
 \begin{align}
 \bfq(x,t) = \sum_i^m\chi_{\mathcal{C}_i}(x)\left[Q_i(t)+(x-x_i)\sigma_i(t)\right],
 \end{align}
where $x_i$ is the center of the cell $\mathcal{C}_i$,  $x_\infty$ the right interface
of $\mathcal{C}_m$, and $\sigma_i(t)$ is a vector of slopes.

We can write \eqref{eq: nonconservative integral} at
a point $x\in \mathcal{C}_j$ ($x<x_j$ for simplicity) as
\begin{align}
    \partial_t\bfq + \partial_x\bff(\bfq)&= \int_x^{x_j} \bfM\partial_z \bfq\,dz
    +\sum_{i=j}^{m-1} \int_{x_i}^{x_{i+1/2}^-} \bfM\partial_z \bfq\,dz\\
    & \quad+ \int_{x_{i+1/2}^+}^{x_{i+1}}  \bfM\partial_z \bfq\,dz
    +\int_{\{x_{i+1/2}\}} \bfM\partial_z \bfq\,dz,\\
    &= \int_x^{x_j} \bfM \sigma_i\,dz
    +\sum_{i=j}^{m-1} \int_{x_i}^{x_{i+1/2}^-} \bfM \sigma_i\,dz\\
    &\quad 
    + \int_{x_{i+1/2}^+}^{x_{i+1}}  \bfM \sigma_{i+1}\,dz
    +\int_{\{x_{i+1/2}\}} \bfM\partial_z \bfq\,dz, \label{eq: split integrals}
\end{align}
where ${x_{i+1/2}^-}$ $x_{i+1/2}^+$ represent the left and right limits of the interface point $x_{i+1/2}$.
Simpson's rule is used to approximate the interior integrals in Equation \eqref{eq: split integrals}.
Since $\bfq$ is discontinuous at $x_{i+1/2}$, we interpret
the nonconservative product as a Borel measure $\mu$  such that
\begin{align}
\label{eq: Borel measure nonconservative}
    \int_{\{x_{i+1/2}\}} \bfM\partial_z \bfq\,dz := \mu({\{x_{i+1/2}\}}) = \int_0^1 \bfM(\phi(s;\bfq^-,\bfq^+)) \frac{\partial \phi}{\partial s}(s;\bfq^-,\bfq^+) ds,
\end{align}
where $\phi$ is a Lipschitz continuous function connecting $ \bfq^{-}(t)=\lim_{x\to x_{i+1/2}^-} \bfq(x,t)$ with 
$ \bfq^+(t)=\lim_{x\to x_{i+1/2}^+}\bfq(x,t)$ \cite{le_floch_shock_1989}. 
Here we choose $\phi$ to be simply a straight line $\phi(s)=s \bfq^+ +(1-s)\bfq^-$, and approximate the 
integral in \eqref{eq: Borel measure nonconservative} with the midpoint rule.

\section{Entropy stability of the slowing-down operator}
\label{sec: Appendix stability slowing down}
One can show that, if only outgoing waves are modified, the slowing-down operator does not destroy 
the entropy-dissipative nature of linear hyperbolic systems.
Indeed, we can consider a general strictly hyperbolic linear system of $N$ PDEs
\begin{equation}
    \partial_t\bfq + \bfA \partial_x\bfq = \bf0,
\end{equation}
defined over the real line such that $\lim_{x\to\pm \infty}\bfq(t,x)=\bf0$.
Let $\bfH$ be an SPD matrix such that $\bfH\bfA$ is symmetric.
Then $\bfv= \frac12 \bfq^T \bfH \bfq$ is a convex entropy function
for the system, 
and $\partial_t \int_\real \bfv(\cdot,t) \,dx=0$.
Let us analyze the effect of the slowing-down operator applied over
the sponge layer $[x_s,x_\infty]$.
For instance, we slow-down the right-going waves, i.e., we solve 
\begin{equation}
\label{eq: appendix slowing}
    \partial_t\bfq + \bfS(x) \partial_x\bfq = \bf0,
\end{equation}
where $\bfS(x)=\bfR \boldsymbol{\Lambda}_\bfS(x) \bfR^{-1}$ such that
$\bfA = \bfR \boldsymbol{\Lambda} \bfR^{-1}$, and the main entries of the diagonal
matrix
$\boldsymbol{\Lambda_S}(x)$ are given by
\begin{equation}
    \tilde{\lambda}_i(x) = \begin{cases}
        \lambda_i, \quad &\text{if } \lambda_i \leq 0,\\
        \Gamma_s(x)\lambda_i, \quad &\text{if } \lambda_i > 0.
    \end{cases}
\end{equation}
The scalar smooth slowing-down function $s(x)$ strictly decreases from $1$ to $0$
inside the sponge layer $[x_s,x_\infty]$ and is constant everywhere else.

Since $\bfA$ is strictly hyperbolic, one can show that $\bfR^{T}\bfH\bfR$ is 
diagonal and positive definite.
Using this fact, we can prove that $\bfH\bfS(x)$ is symmetric and
$\bfH\bfS'(x)$ is symmetric and negative semi-definite for all $x\in \real$, since ${\tilde{\lambda}_i}'(x) \leq 0$.
Therefore, multiplying \eqref{eq: appendix slowing} by $\bfq^T\bfH$ and integrating in space, we get
\begin{align}
    \partial_t \int_\real \bfv(\cdot,t) \,dx &= -\int_\real \bfq^T \bfH \bfS(x) \partial_x\bfq \,dx \\
    &=\frac12 \int_\real \bfq^T \bfH\bfS'(x) \bfq - \partial_x\left(\bfq^T \bfH \bfS(x)\bfq\right)\,dx \\
    &= \frac12 \int_{x_s}^{x_\infty} \bfq^T \bfH\bfS'(x) \bfq\, dx\leq 0.
\end{align}

The same cannot be said in the nonlinear setting,
even for a scalar PDE.
For instance, we can apply the slowing-down operator, described above, to Burgers' equation
\begin{equation}
    \partial_t u + \frac12\partial_x u^2 =0,
\end{equation}
for which we know that $\partial_t ||u||_{L^2(\real)}\leq 0$ holds.
The modified equation reads
\begin{equation}
    \partial_t u + s(x)\frac12\partial_x u^2 =0,
\end{equation}
and after some manipulations, it yields
\begin{align}
    \frac12 \partial_t \int_\real u^2 \, dx &= -\frac13\int_\real s(x)\partial_x u^3 \,dx
    =-\frac13\int_\real \partial_x[s(x)u^3]-s'(x)u^3 \,dx \\
    &=\frac13\int_{x_s}^{x_\infty} s'(x) u^3 \, dx.
\end{align}
If $u<0$ inside the sponge layer, the last term (thus $\partial_t ||u||_{L^2(\real)}$)
becomes strictly greater than zero and the entropy-dissipative nature of the
original system is lost.
This can be attributed to the fact the the sign of the wave speed in Burgers' equation
changes with the solution itself, and thus the slowing-down operator
cannot be defined in a consistent manner,
although a practical solution is to locally modify $s(x)$ such that $s(x)=1$ if $u<0$.

\begin{figure}[]
        \centering
        \includegraphics[width=.9\textwidth]{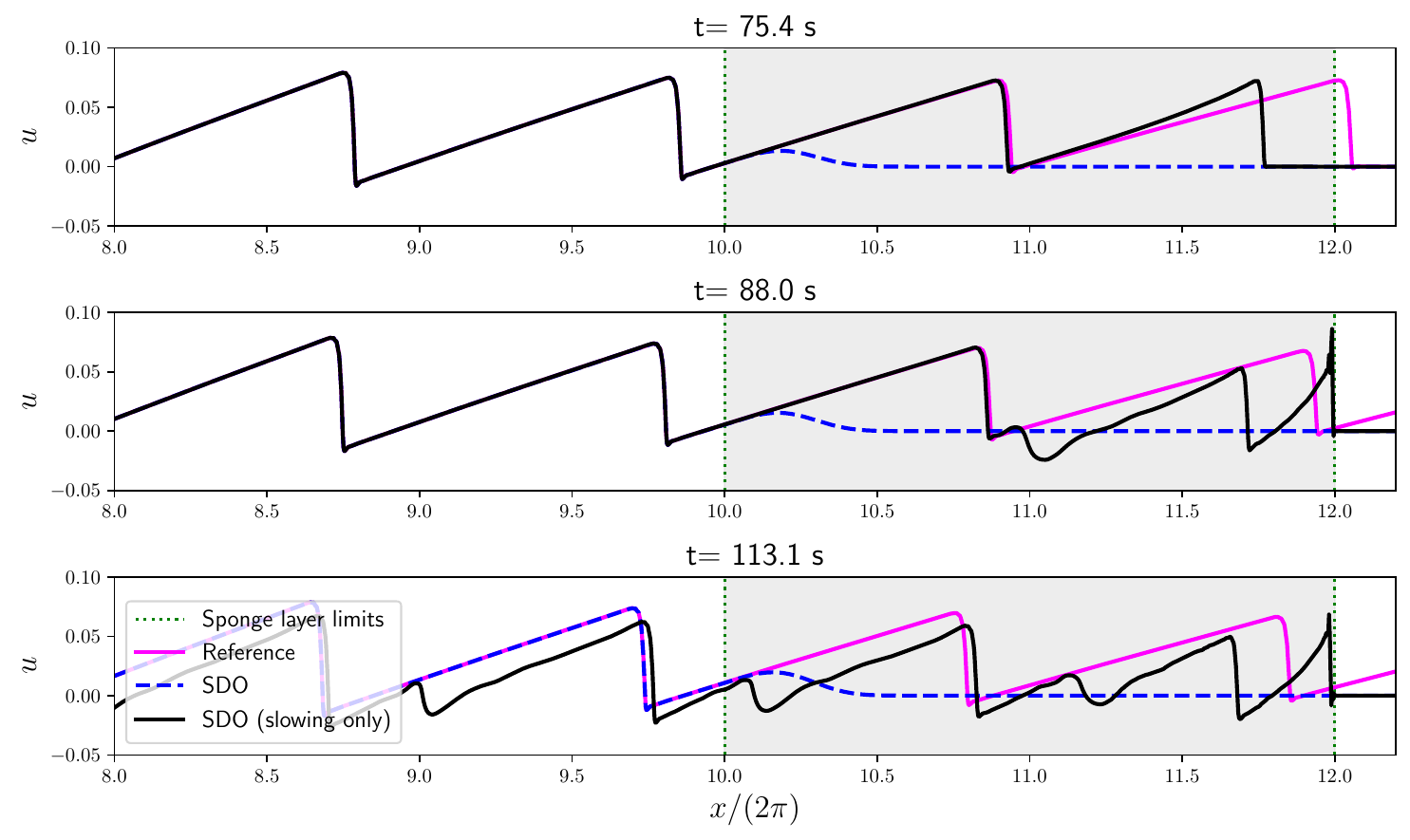}
        \caption[width=.9\textwidth]{Effect of the slowing-down operator 
        for the nonlinear system \eqref{eq:Lagrangian Euler 1D}.
        \label{fig: Appendix entropy 1}}
\end{figure}

\begin{figure}[]
        \centering
        \includegraphics[width=.9\textwidth]{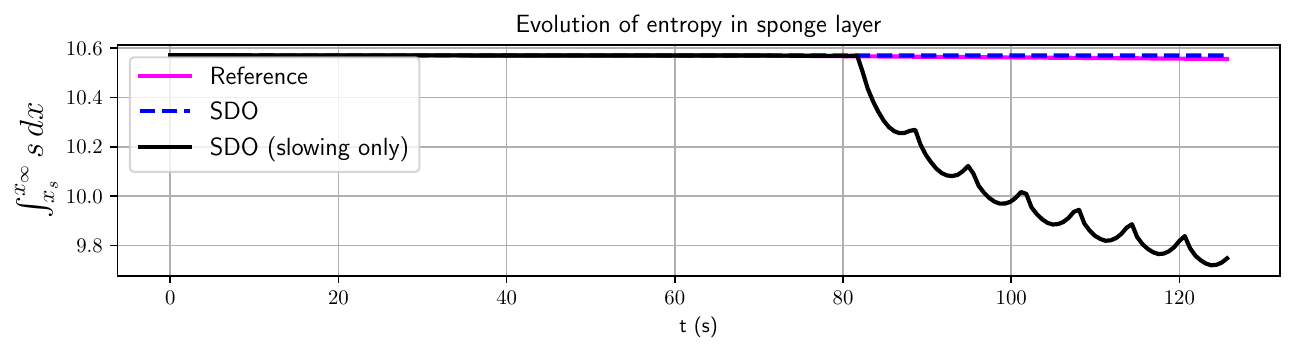}
        \caption[width=.9\textwidth]{Time evolution of the entropy function \eqref{eq:entropy function} integrated over the sponge layer
        $[x_s,x_\infty]$
        \label{fig: Appendix entropy 2}}
\end{figure}

While this is not the case for the nonlinear Euler equations \eqref{eq:Lagrangian Euler 1D},
where the sign of the wave speeds is independent of the solution, 
numerical experiments suggest that the slowing-down operator
can lead to small amounts of entropy production inside the sponge layer.
To illustrate this, we consider the piston problem setup described in Section \ref{sec:Numerics}.
We use $N=250$ cells per wavelength $L$, and a sponge layer of length $\omega=2L$.
To isolate the effect of the slowing-down operator,
we run a simulation with both slowing-down and damping operators,
and one with only the slowing-down operator only.
In Figure \ref{fig: Appendix entropy 1} we can see that, in the absence of a damping operator,
slowing-down leads to wave compression inside the sponge layer.
While the shock-capturing nature of the scheme is able to stabilize the solution,
large spurious oscillations propagate back into the domain.
In Figure \ref{fig: Appendix entropy 2} we can see the evolution of the entropy
function (see e.g. \cite[Sec. 6.4]{whitham_linear_1974})
\begin{align}
    \label{eq:entropy function}
    s(\bfq) = \frac{1}{1-\gamma} \log\left(p V^\gamma\right),
\end{align}
integrated over the sponge layer $[x_s,x_\infty]$.
When only the slowing-down operator is used,
there are periodic increments of entropy inside the sponge layer,
which are controlled by the numerical viscosity
of the discretization.
However, as numerical diffusion vanishes with mesh refinement or
the use of higher-order schemes, this behavior could lead to overall instability.

\end{appendices}
\bibstyle{sn-mathphys-num}
\bibliography{Source_files/ref_Piston_problem}
\end{document}